\documentclass[11pt]{article}
\usepackage{longtable} 
\usepackage{makecell}
\usepackage{array}
\usepackage{hyperref}
\usepackage{amsmath,amsthm,mathtools}
\usepackage{graphicx}
\usepackage{caption}
\usepackage{enumerate}
\usepackage{enumitem}
\setlist[enumerate,1]{label=(\arabic*)}

\usepackage{amsfonts,amssymb,color}

\usepackage[nameinlink,noabbrev]{cleveref}
\crefname{figure}{Fig.}{Figs.}

\newtheorem{theorem}{Theorem}[section]
\newtheorem{lemma}[theorem]{Lemma}

\newtheorem{corollary}[theorem]{Corollary}
\newtheorem{proposition}[theorem]{Proposition}

\theoremstyle{definition}

\newtheorem{remark}[theorem]{Remark}
\newtheorem{problem}[theorem]{Problem}
\newtheorem{example}[theorem]{Example}

\DeclareMathOperator{\gram}{Gram}
\DeclareMathOperator{\supp}{supp}
\DeclareMathOperator{\sspan}{span}

\def\qed{\hfill \rule{4pt}{7pt}}

\begin{document}

\title{\textbf{A Cameron--Goethals--Seidel--Shult theorem\\ for marked graphs}}
\author{ Chaohui Chen$^{a}$, Muhuo Liu$^{b}$,\,   Haiying Shan$^{a}$,  Zoran Stani\'c$^{c}$ \\\small $^{a}$ School of Mathematical Sciences,   Tongji University,   Shanghai,    China\\\small $^{b}$ Department of Mathematics,  South China  Agricultural  University,  Guangzhou, 510642, China \\ \small $^{c}$ Faculty of Mathematics, University of Belgrade, Studentski trg 16, 	11 000 Belgrade, Serbia}

\date{} \maketitle
\noindent
{\bf Abstract.} The classical Cameron--Goethals--Seidel--Shult theorem characterizes graphs
with least eigenvalue at least $-2$ in terms of generalized line graphs,
exceptional graphs, and root systems. In this paper, we establish a marked graph extension of this theorem.
For a graph $G$ on $n$ vertices and a binary vector
$\mathbf{b}\in\{0,1\}^n$, we completely characterize all marked graphs
$\mathbf{G}=(G,\mathbf{b})$ such that
$
     \lambda_{\min}\big(A(G)+\mathbf{b}\mathbf{b}^\top-J_n\big)\ge -2,
$
where $J_n$ denotes the all-ones matrix of order $n$.

Using representations by finite root systems, we show that every irreducible marked graph satisfying the above inequality is represented by either the root system $B_n$ or the exceptional root system $E_8$. This leads to a complete structural characterization in terms of irreducible marked graphs, namely induced subgraphs of extended generalized line graphs, the marked graph $(2K_1,\mathbf{0})$, and exceptional graphs with all vertices marked.

In the special case $\mathbf{b}=\mathbf{j}$, our characterization reduces to the original Cameron--Goethals--Seidel--Shult theorem. The root system $B_n$ and extended generalized line graphs play the roles occupied by $D_n$ and generalized line graphs in the classical theory.
We further apply the characterization
to graphs whose second largest eigenvalue is at most $1$, showing that every $-2$-markable graph satisfies this spectral condition.
 \bigskip


\let\thefootnote\relax\footnotetext{This
	work is    partially supported by the Natural Science Foundation of Guangdong Province (No.~2024A1515011899),
	the National Natural Science Foundation of China (No.~12271182), and the Serbian Ministry of Education, Science and Technological Development via the University of Belgrade.} 

\let\thefootnote\relax\footnotetext{
	\textit{E-mail addresses}: \texttt{cch956833210@163.com} (C. Chen),  \texttt{liumuhuo@163.com} (M. Liu),   \texttt{shan\_haiying@tongji.edu.cn} (H. Shan,  Corresponding author),  \texttt{zstanic@matf.bg.ac.rs} (Z. Stani\'c)}

\let\thefootnote\relax\footnotetext{\textit{ORCID IDs}: \texttt{0000-0002-4586-3120} (C. Chen),  \texttt{0000-0001-8217-2452} (M. Liu),   \texttt{0000-0002-2749-9476} (H. Shan), \texttt{0000-0002-4949-4203} (Z. Stani\'c)}

 \section{Introduction}
All graphs considered are finite and undirected. Unless stated otherwise, they are assumed to be simple, that is, without loops or multiple edges. The \emph{eigenvalues} of a graph $G$ of order $n$ are the eigenvalues of its adjacency matrix $A(G)$, indexed in non-increasing order as
$
\lambda_1 \ge \lambda_2 \ge \cdots \ge \lambda_n.
$
When the order is irrelevant or clear from the context, we write $\lambda_{\min}$ for the least eigenvalue. The same notation is used for arbitrary real symmetric matrices.  We write $I$, $J$, $O$, $\mathbf{j}$, and $\mathbf{0}$ for the identity matrix, the all-ones matrix, the all-zero matrix, the all-ones vector, and the all-zero vector, respectively, with the size specified by a subscript when necessary.

The following problem is a cornerstone of spectral graph theory and has attracted sustained attention since the 1960s.

\begin{problem}\label{pro0}\cite{Hof1968}
Characterize all  graphs $G$ that satisfy
$
 \lambda_{\min}(G) \geq -2.
$
\end{problem}

In the context of Problem~\ref{pro0}, connected graphs with least eigenvalue at least $-2$ are traditionally divided into two disjoint families: generalized line graphs and exceptional graphs; their definitions are recalled in the next section. A major breakthrough was achieved in 1976 when Cameron, Goethals, Seidel, and Shult~\cite{Cam1976} established a structural classification of such graphs using the theory of root systems. Their theorem shows that generalized line graphs are precisely the graphs representable in the root system $D_n$, whereas exceptional graphs are those that are not representable in any $D_n$ and arise from the root system $E_8$. Moreover, every exceptional graph has at most $36$ vertices and maximum degree at most $28$. A complete classification of exceptional graphs was obtained in the late 1990s independently by several authors. Further details and historical remarks can be found in \cite{Cvet2004}.

Motivated by this theory, it is natural to investigate a mild yet meaningful perturbation of the adjacency matrix. Let $\mathbf{b}\in\{0,1\}^n$ be a binary vector and consider the symmetric matrix $A(G)+\mathbf{b}\mathbf{b}^{\top}-J_n$. This leads to the following generalization of Problem~\ref{pro0}.

\begin{problem}\label{pro1}
Characterize all pairs \((G, \mathbf{b})\), where \(G\) is a graph on \(n\) vertices and \(\mathbf{b} \in \{0,1\}^n\), such that
\[
 \lambda_{\min}(A(G) + \mathbf{b}\mathbf{b}^\top - J_n) \geq -2.
\]
\end{problem}

The aim of this paper is to study Problem~\ref{pro1} from a representation-theoretic perspective. Observe that when $\mathbf{b}$ is the all-ones vector $\mathbf{j}$, the matrix $A(G)+\mathbf{b}\mathbf{b}^{\top}-J_n$ reduces to $A(G)$, and thus Problem~\ref{pro1} reduces to Problem~\ref{pro0}. In the classical setting, the root systems $D_n$ and $E_8$ arising in the structural classification consist of roots of a single length. In contrast, the introduction of the binary vector $\mathbf{b}$ naturally leads to the root system $B_n$, which contains roots of two distinct lengths and contains $D_n$ as a subsystem. Consequently, the representing vectors no longer have a common length but instead form a configuration involving two different lengths.

In this paper, we establish a marked graph extension of the Cameron--Goethals--Seidel--Shult theorem and thereby obtain a complete solution to Problem~\ref{pro1}. The main results are stated in Theorems~\ref{thm:-2rep} and~\ref{thm2}. These theorems extend the characterization, summarized in \cite{Cvet2004}, of graphs arising in Problem~\ref{pro0}. More generally, Theorem~\ref{thm2} is closely related to results on graphs with bounded least or largest eigenvalues. This topic has been extensively studied; consequently, we omit further discussion and instead refer the reader to the monographs \cite{Cvet2004,ifge}. For recent developments and generalizations, see \cite{Cvet2015,Koo2021} and \cite{Bel2016,Grea2015,Row2022}, respectively.

Section~\ref{sec2} introduces the necessary definitions and terminology and states the main results of the paper. Section~\ref{sec3} may be viewed as a continuation of Section~\ref{sec2}, focusing on root systems. Although much of the material in this section is standard, it is included to establish notation and recall results needed later. In Section~\ref{sec4} we develop a representation theory for marked graphs based on root systems. Section~\ref{sec5} is devoted to the proofs of Theorems~\ref{thm:-2rep} and~\ref{thm2}. Finally, Section~\ref{sec6} applies the characterization of Theorem~\ref{thm2} to the study of $-2$-markable graphs and their relationship with graphs satisfying $\lambda_2 \le 1$.

 \section{Main  result}\label{sec2}
 
 Real vectors are denoted by bold Latin letters or by Greek letters in upright type, depending on the context. In particular, $\{\mathbf{e}_i\}_{i=1}^{n}$  denotes the canonical orthonormal basis of $\mathbb{R}^n$. As usual, $K_n$ and $P_n$  denote the  complete graph and the path of order $n$, respectively. Also,  $G_1\cup G_2$ denotes the disjoint union of  $G_1$ and $G_2$, and $G_1\vee G_2$ denotes the join of $G_1$ and $G_2$ (by convention, the vertices of $G_1$ are ordered before those of $G_2$ in the overall labeling). 

Now, for any graph \(G\) on \(n\) vertices and any binary vector \(\mathbf{b} \in \{0,1\}^n\), we label the vertices of \(G\) consistently with the indexing of \(\mathbf{b}\). We denote by \(V_{\mathbf{b}}(G)\) the \textit{marked vertex set} of \(G\), which consists of all vertices whose corresponding entry in \(\mathbf{b}\) is equal to 1.

For convenience, we write \(\mathbf{G} = (G, \mathbf{b})\) to denote the \textit{marked graph} with underlying graph \(G\) and marking vector \(\mathbf{b}\). Note that the same underlying graph \(G\) can give rise to different marked graphs by choosing different marking vectors \(\mathbf{b}\). For drawing purposes, marked vertices of a  marked graph $\mathbf{G}=(G,\mathbf{b})$ are conventionally coloured black, while unmarked vertices are coloured white.

For a marked graph \(\mathbf{G} = (G, \mathbf{b})\), we define the \textit{marked adjacency matrix} \(\hat{A}(\mathbf{G})\) of \(\mathbf{G}\) by
\[
\hat{A}(\mathbf{G}) = A(G) + \mathbf{b}\mathbf{b}^\top - J_n.
\]

Two marked graphs \(\mathbf{G}_1 = (G_1, \mathbf{b}_1)\) and \(\mathbf{G}_2 = (G_2, \mathbf{b}_2)\) are called \textit{isomorphic}, written \(\mathbf{G}_1 \cong \mathbf{G}_2\), if there exists a graph isomorphism \(\varphi \colon V(G_1) \to V(G_2)\) that preserves the marked vertex sets, that is, \(\varphi(V_{\mathbf{b}_1}(G_1)) = V_{\mathbf{b}_2}(G_2)\).
A   marked graph \( \mathbf{G}'=(G', \mathbf{b}') \) is called an \textit{induced subgraph} of \( \mathbf{G}=(G, \mathbf{b}) \) if \( G' \) is an induced subgraph of \( G \) and \( V_{\mathbf{b}'}(G')=V_{\mathbf{b}}(G)\cap V(G') \).

Let $\mathbf G_1=(G_1,\mathbf b_1)$ and $\mathbf G_2=(G_2,\mathbf b_2)$ be
two vertex-disjoint marked graphs, and let
\[
    \mathbf b=
    \begin{pmatrix}
        \mathbf b_1\\
        \mathbf b_2
    \end{pmatrix}.
\]
The \textit{marked union} $\mathbf G_1\widehat{\cup}\mathbf G_2$ is the
marked graph obtained from
$ \left(G_1\vee G_2,\mathbf b\right)$
by deleting all edges between $V_{\mathbf b_1}(G_1)$ and
$V_{\mathbf b_2}(G_2)$. Thus, if
$\mathbf b_1=\mathbf j$ and $\mathbf b_2=\mathbf j$, then
$
\mathbf G_1\widehat{\cup}\mathbf G_2=(G_1\cup G_2,\mathbf b);
$
if $\mathbf b_1=\mathbf 0$ or $\mathbf b_2=\mathbf 0$, then
$
\mathbf G_1\widehat{\cup}\mathbf G_2=(G_1\vee G_2,\mathbf b).
$

Dually, the \textit{marked join} $\mathbf G_1\widehat{\vee}\mathbf G_2$ is
the marked graph obtained from
$
    \left(G_1\cup G_2,\mathbf b\right)
$
by adding all edges between $V_{\mathbf b_1}(G_1)$ and
$V_{\mathbf b_2}(G_2)$. Thus, if $\mathbf b_1=\mathbf j$ and
$\mathbf b_2=\mathbf j$, then
$
\mathbf G_1\widehat{\vee}\mathbf G_2=(G_1\vee G_2,\mathbf b);
$
if $\mathbf b_1=\mathbf 0$ or $\mathbf b_2=\mathbf 0$, then
$
\mathbf G_1\widehat{\vee}\mathbf G_2=(G_1\cup G_2,\mathbf b).
$

\begin{figure}[htbp]
  \centering
  \includegraphics[page=1,scale=0.99]{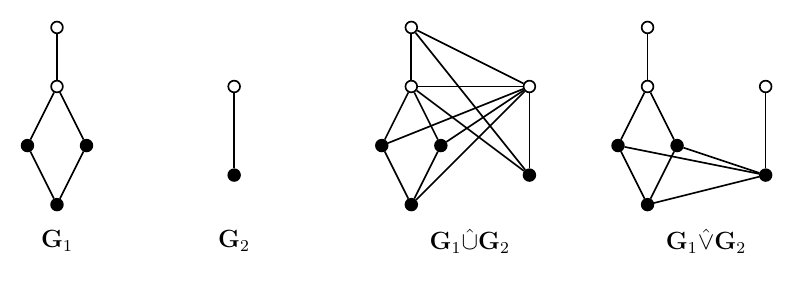}
  \caption{An illustration of the marked union and marked join operations.}
  \label{fig:marked-union-join}
\end{figure}

Figure~\ref{fig:marked-union-join} illustrates the marked union and marked
join operations for two marked graphs.

More generally, for a finite family of pairwise vertex-disjoint marked graphs
$\{\mathbf G_i=(G_i,\mathbf b_i)\}_{i=1}^r$, we define
\[
    \widehat{\cup}_{i=1}^1 \mathbf G_i:=\mathbf G_1,\quad
    \widehat{\cup}_{i=1}^r \mathbf G_i
    :=
    \left(\widehat{\cup}_{i=1}^{r-1}\mathbf G_i\right)
    \widehat{\cup}\mathbf G_r
    \quad (r\ge 2),
\]
and similarly
\[
    \widehat{\vee}_{i=1}^1 \mathbf G_i:=\mathbf G_1,\quad
    \widehat{\vee}_{i=1}^r \mathbf G_i
    :=
    \left(\widehat{\vee}_{i=1}^{r-1}\mathbf G_i\right)
    \widehat{\vee}\mathbf G_r
    \quad (r\ge 2).
\]
Consider two marked graphs
$\mathbf G_1=(G_1,\mathbf b_1)$ and $\mathbf G_2=(G_2,\mathbf b_2)$. By the
definition of marked union, the off-diagonal blocks of the adjacency matrix
of the underlying graph of $\mathbf G_1\widehat{\cup}\mathbf G_2$ are
$
    J-\mathbf b_1\mathbf b_2^\top
    \quad\text{and}\quad
    J-\mathbf b_2\mathbf b_1^\top,
$
with all-ones matrices 
 of appropriate sizes. Hence,
\[
\begin{aligned}
\hat{A}(\mathbf G_1\widehat{\cup}\mathbf G_2)
&=
\begin{pmatrix}
A(G_1) & J-\mathbf b_1\mathbf b_2^\top\\
J-\mathbf b_2\mathbf b_1^\top & A(G_2)
\end{pmatrix}
+
\begin{pmatrix}
\mathbf b_1\\
\mathbf b_2
\end{pmatrix}
\begin{pmatrix}
\mathbf b_1\\
\mathbf b_2
\end{pmatrix}^{\top}
-J  \\
&=
\begin{pmatrix}
A(G_1)+\mathbf b_1\mathbf b_1^\top-J & O\\
O & A(G_2)+\mathbf b_2\mathbf b_2^\top-J
\end{pmatrix}  \\
&=
\begin{pmatrix}
\hat{A}(\mathbf G_1) & O\\
O & \hat{A}(\mathbf G_2)
\end{pmatrix}.
\end{aligned}
\]

More generally, for a finite family of pairwise vertex-disjoint marked graphs
$\{\mathbf G_i\}_{i=1}^r$, we obtain
\begin{equation}\label{eqq}
\hat{A}\left(\widehat{\cup}_{i=1}^r \mathbf G_i\right)
=
\begin{pmatrix}
\hat{A}(\mathbf G_1) & O & \cdots & O\\
O & \hat{A}(\mathbf G_2) & \cdots & O\\
\vdots & \vdots & \ddots & \vdots\\
O & O & \cdots & \hat{A}(\mathbf G_r)
\end{pmatrix}.
\end{equation}

A   marked graph $\mathbf{G}$ is \textit{irreducible} if $\mathbf{G}$ is not the marked union of any two  marked graphs; equivalently, the marked matrix \(\hat{A}(\mathbf{G})\) is not permutation similar to a direct sum of two non-trivial proper principal submatrices. Accordingly, every   marked graph $\mathbf{G}$ has the unique \textit{irreducible decomposition}
\[\mathbf{G}\cong\widehat{\cup}_{i=1}^{r} \mathbf{G}_i,\]
where $r\ge 1$ and each $\mathbf{G}_i$ is irreducible.

\begin{remark}
For $\mathbf{b}=\mathbf{j}$, we have $\hat{A}(\mathbf G)=A(G)$. Moreover, in the all-marked case, marked union coincides with ordinary disjoint union; hence the irreducibility of $\mathbf G$ is equivalent to the connectedness of $G$.
\end{remark}

At this point, we recall the definitions of line graphs, generalized line graphs, and exceptional graphs. Although these notions are well known, our subsequent arguments require that they should be stated in a precise and consistent form.

The \textit{line graph} $L(H)$ of a multigraph (i.e., graph that admits parallel edges, but with no loops) $H$ is defined as follows. The vertices of $L(H)$ are the edges of $H$, and two vertices of $L(H)$ are adjacent if and only if  the corresponding edges of $H$ have exactly one common vertex in $H$.

\begin{figure}[h!]
  \centering
  \includegraphics[page=2,scale=1]{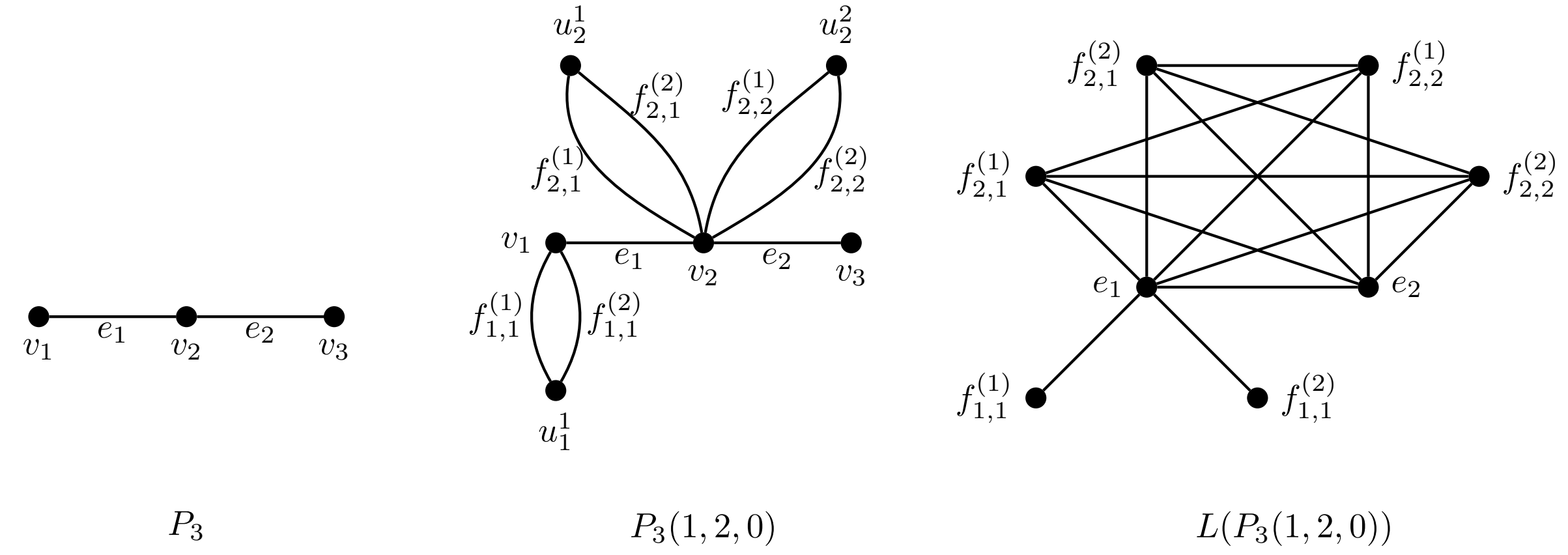}
  \caption{The construction of a generalized line graph.}
  \label{fig:3}
\end{figure}

The following definition of the generalized line graph is adapted from \cite{Cvet2010}. A pendant double edge is called a \textit{petal}, while a \textit{blossom} $B^m$ consists of $m$ petals ($m \geq 0$) attached to a single vertex. The vertices of the blossom distinct from this common vertex are called \textit{petal vertices} (An empty blossom $B^0$ has no petals and reduces to the trivial graph $K_1$). For a graph $G$ with vertex set $V(G) = \{v_1, v_2, \dots, v_n\}$ and non-negative integers $a_1, a_2, \ldots, a_n$, the \textit{$B$-graph} $H \cong G(a_1,a_2, \ldots, a_n)$ is the multigraph obtained by attaching a blossom $B^{a_i}$ to $v_i$ of $G$. To make a clear  distinction between vertices of $G$ and those introduced by blossoms in $H$, we relabel the vertices of $H$ as follows: For each vertex $v_i \in V(G)$, we label the petal vertices attached at $v_i$ as $u_i^1,u_i^2, \ldots, u_i^{a_i}$. The line graph $L(H)$, also denoted by $L(G; a_1, a_2, \ldots, a_n)$, is called a \textit{generalized line graph}. The construction of a generalized line graph is illustrated in Figure~\ref{fig:3}.

For \(H\cong G(a_1, a_2,\dots,a_n)\), let \(u_i^{\,j}\) be the \(j\)-th petal vertex attached to \(v_i\) \((1\le j\le a_i)\).
The two parallel edges between \(v_i\) and \(u_i^{\,j}\) are denoted by  \(f_{i,j}^{(1)}\) and \(f_{i,j}^{(2)}\), respectively.
  Let $U^p$ be the set of all petal vertices in $H$. Then $V(H)=V(G)\cup U^p$. Define the incidence matrix \(M\in\{-1,0,1\}^{\,|V(H)|\times|E(H)|}\) by indexing rows by \(v\in V(H)\) and columns by \(e\in E(H)\) and setting
\[
M_{v,e}=
\begin{cases}
-1,& \text{if }(v,e)=(u_i^{\,j},\,f_{i,j}^{(2)}),\\[2pt]
1,& \text{if } v\in e\ \text{and }(v,e)\ne(u_i^{\,j},\,f_{i,j}^{(2)}),\\[2pt]
0,& \text{if } v\notin e.
\end{cases}
\]
In other words, if \(\{\mathbf{e}_v\}_{v \in V(H)}\) is the standard orthonormal basis of \(\mathbb{R}^{|V(H)|}\), then 
the column vectors of \(M\) are precisely the vectors 
\begin{equation}\label{eqgl}
    \{\mathbf{e}_{v_i}+\mathbf{e}_{v_j}\,:\,(v_i,v_j)\in E(G)\}\cup\{\mathbf{e}_{v_i}\pm\mathbf{e}_{u_i^j}\,:\,u_i^j\in U^p\}.
\end{equation}

For example, the incidence matrix of the $B$-graph $P_3(1,2,0)$ from Figure~\ref{fig:3} is
\[M=\begin{pmatrix}
    1 & 0 & 1 & 1 & 0 & 0 & 0 & 0\\
    1 & 1 & 0 & 0 & 1 & 1 & 1 & 1\\
    0 & 1 & 0 & 0 & 0 & 0 & 0 & 0\\
    0 & 0 & 1 & -1 & 0 & 0 & 0 & 0\\
    0 & 0 & 0 & 0 & 1 & -1 & 0 & 0\\
    0 & 0 & 0 & 0 & 0 & 0 & 1 & -1\\
\end{pmatrix},\]
where the rows are indexed by $v_1,v_2,v_3,u_1^1,u_2^1,u_2^2$ and the columns are indexed by $e_1,e_2,f^{(1)}_{1,1},f^{(2)}_{1,1},$ $f^{(1)}_{2,1},f^{(2)}_{2,1},f^{(1)}_{2,2},f^{(2)}_{2,2}$.

With the incidence matrix $M$ defined above, we have $A(L(H))+2I=M^\top M$. This implies that $ \lambda_{\min}(A(L(H)))\ge -2$. That is, every generalized line graph satisfies $ \lambda_{\min}\ge -2$.
A connected graph $G$ that is not a generalized line graph but satisfies $ \lambda_{\min}(G)\ge -2$ is called an \textit{exceptional graph}. Consequently, every connected graph $G$ with $ \lambda_{\min}(G)\ge -2$ is either a generalized line graph or an exceptional graph.

We now introduce further definitions. Let $G$ be a graph with vertex set $V(G)=\{v_1,v_2,\dots,v_n\}$ and let  $H\cong G(a_1,a_2,\dots,a_n)$ be the $B$-graph associated with $G$. The graph $L^*(H)$ is obtained from the disjoint union of $L(H)$ and the complete graph $K_n$ on the vertex set $\{v_1,v_2,\dots,v_n\}$ by adding edges between a vertex $e\in V(L(H))$ and a vertex $v_i\in V(K_n)$ if and only if the edge $e$ of $H$ is not incident with $v_i$ in $H$.
In other words, put \(m=|E(H)|\), and assume that \( M \) is an incidence matrix of \( H \) with the block form
\[
M = \begin{pmatrix} M_1 \\ M_2 \end{pmatrix},
\]
where the rows of \( M_1\in\{0,1\}^{n\times m} \) correspond to the vertices in \( V(G) \), and the rows of \( M_2 \in\{-1,0,1\}^{|U^p|\times m}\) correspond to the petal vertices. Then the adjacency matrix of $L^*(H)$ has the form
\[
A\big(L^*(H)\big)=
\begin{pmatrix}
A\!\big(L(H)\big) & J_{m\times n}-M_1^\top\\
J_{n\times m}-M_1 & J_{n}-I_{n}
\end{pmatrix}.
\]
Let $\mathbf{L}^*(H)=\left(L^*(H),\begin{pmatrix}
    \mathbf{j}_{m}\\
    \mathbf{0}_{n}
\end{pmatrix}\right)$, that is, $\mathbf{L}^*(H)$ is a  marked graph of $L^*(H)$ with the marked vertex  set   $V(L(H))$. We call $\mathbf{L}^*(H)$ the \textit{extended generalized line graph} of the $B$-graph $H$.

The importance of extended generalized line graphs stems from the fact that they naturally admit representations containing both long and short roots of the root system $B_n$. As will be shown below, they play the same role in the theory of marked graphs as generalized line graphs do in the classical Cameron--Goethals--Seidel--Shult theorem. The construction of $\mathbf{L}^*(H)$ for $H\cong P_3(1,2,0)$ (shown in Figure~\ref{fig:3}) is illustrated in Figure~\ref{fig:4}.

\begin{figure}[h!]
  \centering
  \includegraphics[page=3,scale=0.9]{fig.pdf}
  \caption{The construction of an extended generalized line graph.}
  \label{fig:4}
\end{figure}

By considering the matrix $M'=\begin{pmatrix}
    M_1 & -I_n\\
    M_2 & O_{  |U^p|\times n}
\end{pmatrix}$, we obtain
\begin{equation}\label{eqcl}
    \hat{A}(\mathbf{L}^*(H))+2I=\begin{pmatrix}
A\!\big(L(H)\big)+2I & -M_1^\top\\
-M_1 & I_n
\end{pmatrix}=\begin{pmatrix}
M^\top M & -M_1^\top\\
-M_1 & I_n
\end{pmatrix}=M'^\top M',
\end{equation}
along with
\[
 \lambda_{\min}(\hat{A}(\mathbf{L}^*(H)))\ge -2.
\]

An $n \times n$ matrix $M$ is \textit{represented} by a system of real vectors $\Phi$ if there exists 
an \(n\)-tuple
\(\Sigma=(\alpha_1,\alpha_2,\ldots,\alpha_n)\) of vectors in \(\Phi\) such that its Gram matrix satisfies
\begin{equation}\label{41e}
    \gram(\Sigma) = M + 2I_n.
\end{equation}
A marked graph $\mathbf{G} = (G, \mathbf{b})$ with vertex set $V(G) = \{1, 2, \ldots, n\}$ is said to be \textit{represented} by $\Phi$ if $\hat{A}(\mathbf{G})$ is represented by $\Phi$. In this case, $\Sigma$ is called a \textit{representation} of $\mathbf{G}$.  Since the vectors in $\Sigma$ are always distinct in our applications, we generally write $\Sigma$ as a tuple, but occasionally as a set when convenient. An implicit ordering is assumed when needed.

\begin{remark}\label{rem:repT2}   Let \(H\cong G(a_1, a_2, \ldots, a_n)\) be a $B$-graph with $V(H)=V(G)\cup U^p$ and incidence matrix $M$. Since $A(L(H))+2I=M^\top M$, $A(L(H))$ is represented by the tuple formed by the column vectors of \(M\), whose entries are precisely the vectors in \eqref{eqgl}. Also, by \eqref{eqgl} and~\eqref{eqcl}, the tuple \( \Sigma \) formed by the column vectors of \( M' \), indexed by the columns of \(M'\), is a representation of $\mathbf{L}^*(H)$; its entries are precisely the vectors  
\[
\bigl\{ \mathbf{e}_{v_i} + \mathbf{e}_{v_j} \,:\, (v_i, v_j) \in E(G) \bigr\}\cup \bigl\{ \mathbf{e}_{v_i} \pm \mathbf{e}_{u_i^j} \,:\, u_i^j \in U^p \bigr\}\cup \bigl\{ -\mathbf{e}_{v_i} \,:\, v_i \in V(G) \bigr\}.
\]
In particular, when \(G \cong K_n\),  these entries are precisely the vectors  
\[
\bigl\{ \mathbf{e}_{v_i} + \mathbf{e}_{v_j} \,:\, v_i\ne v_j \in V(G)\bigr\} \cup \bigl\{ \mathbf{e}_{v_i} \pm \mathbf{e}_{u_i^j} \,:\, u_i^j \in U^p \bigr\} \cup \bigl\{ -\mathbf{e}_{v_i} \,:\, v_i \in V(G) \bigr\}.
\]
\end{remark}

The Cameron--Goethals--Seidel--Shult theorem  is stated as follows, where the explicit descriptions of the root systems $D_n$ and $E_8$ are given in the subsequent Theorem~\ref{rootthm}.
\begin{theorem}[\cite{Cam1976}]\label{thm:Cam}
    Let $G$ be a connected graph. Then the following statements hold:
    \begin{enumerate}
        \item[(i)] $ \lambda_{\min}(A(G))\ge -2$ if and only if $A(G)$ is represented by some $D_n$ or $E_8$;
        \item[(ii)] $A(G)$ is represented by some $D_n$ if and only if $G$ is a generalized line graph;
        \item[(iii)] $A(G)$ is represented by $E_8$ and is not represented by any $D_n$ if and only if $G$ is an exceptional graph. 
    \end{enumerate}
\end{theorem}

The following theorem may be viewed as a marked-graph analogue of Theorem~\ref{thm:Cam}, with the root system $B_n$ replacing $D_n$ and extended generalized line graphs replacing generalized line graphs. Again, we refer to Theorem~\ref{rootthm} for the explicit descriptions of the root systems.

\begin{theorem}\label{thm:-2rep}
    Let $\mathbf{G}=(G,\mathbf{b})$ be an irreducible  marked graph. Then the following statements hold:
    \begin{enumerate}
        \item[(i)] $ \lambda_{\min}(\hat{A}(\mathbf{G}))\ge -2$ if and only if $\mathbf{G}$ is represented by some $B_n$ or by $E_8$;
        \item[(ii)] $\mathbf{G}$ is represented by some $B_n$ if and only if $\mathbf{G}$ is an induced subgraph of an extended generalized line graph or $(2K_1,\mathbf{0})$;
        \item[(iii)] $\mathbf{G}$ is represented by $E_8$ and is not represented by any $B_n$ if and only if $\mathbf{b}=\mathbf{j}$ and $G$ is an exceptional graph. 
    \end{enumerate}
\end{theorem}

By removing the irreducibility assumption, the following theorem extends 
to a characterization of all marked graphs with 
$ \lambda_{\min}(\hat{A}(\mathbf{G}))\ge -2$, expressed via the irreducible 
decomposition introduced in Section~\ref{sec2}.
\begin{theorem}\label{thm2} Let $\mathbf{G}=(G,\mathbf{b})$ be a marked graph, and let
\[
    \mathbf{G} = \widehat{\cup}_{i=1}^{r} \mathbf{G}_i
\]
be its irreducible decomposition. Then $ \lambda_{\min}(\hat{A}(\mathbf{G}))\ge -2$ if and only if, for every \(1\le i\le r\), \(\mathbf{G}_i=(G_i,\mathbf{b}_i)\) is one of the following types:
	\begin{itemize}[leftmargin=1.5cm]
		\item[(T1)] $(G_i, \mathbf{b}_i)$ is an induced subgraph of an extended generalized line graph;
		\item[(T2)] \( (G_i, \mathbf{b}_i)=(2K_1, \mathbf{0})\);
		\item[(T3)] $G_i$ is an exceptional graph and $\mathbf{b}_i=\mathbf{j}$.
	\end{itemize}
\end{theorem}

As noted above,  a   marked graph $\mathbf{G}$ satisfies $ \lambda_{\min}\bigl(\hat{A}(\mathbf{G})\bigr) \ge -2$ if and only if the matrix $\hat{A}(\mathbf{G}) + 2I$ is positive semidefinite. In analogy with classical spectral graph theory, the study of graphs satisfying $ \lambda_{\min}(G)\ge -2$ is intrinsically linked to the theory of finite root systems. The necessary background on root systems is recalled in the next section.

\section{Background on root systems}\label{sec3}
We write $\langle \cdot, \cdot \rangle$ for the inner product in a (finite-dimensional) Euclidean space and $\|\cdot\|$ for the corresponding norm. For $r \in \mathbb{R}$ and a set of vectors $S$, we denote by $rS$ the scaled set $\{ r\alpha \,:\, \alpha \in S \}$.

A root system is a fundamental algebraic structure that arises in the representation theory of Lie algebras and Lie groups, consisting of a finite set of non-zero vectors in a Euclidean space satisfying specific symmetry and integrality conditions~\cite{Hal2015,Hum1972}. Precisely, let $\Phi$ be a set of non-zero vectors in a Euclidean space $E$. Then $\Phi$ is called a \textit{root system} in $E$ if the following assumptions hold:
\begin{enumerate}[leftmargin=1.5cm]
    \item[(R1)] $\Phi$ is finite and spans $E$.
    \item[(R2)] If $\alpha \in \Phi$, then the only multiples of $\alpha$ in $\Phi$ are $\pm \alpha$.
    \item[(R3)] If $\alpha,\beta \in \Phi$, then the reflection $\sigma_\alpha(\beta)\coloneqq \beta-2\frac{\langle\alpha,\beta\rangle}{\langle\alpha,\alpha\rangle}\alpha\in \Phi$.
    \item[(R4)] If $\alpha, \beta \in \Phi$, then $2\frac{\langle\alpha,\beta\rangle}{\langle\alpha,\alpha\rangle} \in \mathbb{Z}$.
\end{enumerate}

The root systems $\Phi$, $\Phi'$ in  Euclidean spaces $E$ and $E'$, respectively, are \textit{isomorphic} if there exists an invertible linear transformation $\phi \colon E \to E'$ such that $\phi(\Phi)=\Phi'$ and for all $\alpha,\beta \in \Phi$, we have
    \begin{equation*}\label{eqiso}
        \frac{\langle\phi(\alpha),\phi(\beta)\rangle}{\langle\phi(\alpha),\phi(\alpha)\rangle}=\frac{\langle\alpha,\beta\rangle}{\langle\alpha,\alpha\rangle}.
    \end{equation*}
A map $\phi$ with this property is called an \textit{isomorphism}.
The elements of a root system $\Phi$ are called \textit{roots}, and $l(\Phi)=\{\|\alpha\|\,:\,\alpha\in\Phi\}$ denotes the set of all possible root lengths in $\Phi$.
\begin{remark}\label{remark1}
    Note that the isomorphism $\phi$ is not required to preserve inner products or root lengths. For instance, $\Phi$ is isomorphic to $c\Phi$ for any $c \neq 0$. Nevertheless, $\phi$ does preserve reflections across roots: For all $\alpha, \beta \in \Phi$,
    \[
    \phi\bigl(\sigma_\alpha(\beta)\bigr) = \sigma_{\phi(\alpha)}\bigl(\phi(\beta)\bigr),
    \]
    where $\sigma_\alpha(\beta)$ is defined in (R3).
\end{remark}

A set of non-zero vectors in a Euclidean space $E$ is \textit{irreducible} if it cannot be partitioned into  two subsets such that each vector in one set is orthogonal to every vector in the other.

We now present two lemmas and a characterizing theorem.

\begin{lemma}[\cite{Hum1972}]\label{lem:rootlength}
   \, If $\Phi$ is an irreducible root system, then at most two root lengths occur in $\Phi$, i.e., $|l(\Phi)|\le 2$.
\end{lemma}
\begin{lemma}[\cite{Hum1972}]\label{lemire}
   \, Every root system $\Phi$ in a Euclidean space $E$ admits a unique decomposition $\Phi=\bigcup_{i=1}^{r}\Phi_i$, where $r\ge 1$, $E=\bigoplus_{i=1}^{r}E_i$ is an orthogonal direct sum,  and each $\Phi_i$ is an irreducible root system in the subspace $E_i$.
\end{lemma}

The announced characterization of root systems, well known and appearing in many sources, is as follows.

\begin{theorem}[\cite{Hal2015,Hum1972}]\label{rootthm} Every irreducible root system is isomorphic to exactly one of the following:
\begin{enumerate}
    \item[(i)] $A_n=\{\mathbf{e}_i-\mathbf{e}_j\,:\,i\ne j,\,\, 1\le i,j\le n+1\}$, $n \geq 1$, $\mathbf{e}_i\in\mathbb{R}^{n+1}$;
    \item[(ii)] $B_n=\{\pm\mathbf{e}_i\pm\mathbf{e}_j\,:\,1\le i<j\le n\}\cup\{\pm\mathbf{e}_i\,:\,1\le i\le n\}$, $n \geq 2$, $\mathbf{e}_i\in\mathbb{R}^{n}$;
    \item[(iii)] $C_n=\{\pm\mathbf{e}_i\pm\mathbf{e}_j\,:\,1\le i<j\le n\}\cup\{\pm2\mathbf{e}_i\,:\,1\le i\le n\}$, $n \geq 2$, $\mathbf{e}_i\in\mathbb{R}^{n}$;
    \item[(iv)] $D_n=\{\pm\mathbf{e}_i\pm\mathbf{e}_j\,:\,1\le i<j\le n\}$, $n \geq 3$, $\mathbf{e}_i\in\mathbb{R}^{n}$;
    \item[(v)] $E_8=D_8\cup\big\{\tfrac{1}{2}\sum_{i=1}^{8}r_i\mathbf{e}_i\,:\,r_i=\pm 1,\prod_{i=1}^8r_i=1\big\}$, along with two proper subsets $E_6,E_7$ of $E_8$;
    \item[(vi)] $F_4=B_4\cup\{\tfrac{1}{2}(\pm\mathbf{e}_1\pm\mathbf{e}_2\pm\mathbf{e}_3\pm\mathbf{e}_4)\}$;
    \item[(vii)] $G_2=A_2\cup\{\pm(2\mathbf{e}_1-\mathbf{e}_2-\mathbf{e}_3),\pm(2\mathbf{e}_2-\mathbf{e}_1-\mathbf{e}_3),\pm(2\mathbf{e}_3-\mathbf{e}_1-\mathbf{e}_2)\}$.
\end{enumerate}
\end{theorem}
The restrictions on $n$ are to avoid the case of $D_2$, which is not irreducible. For explicit set-theoretic descriptions of $E_6$ and $E_7$, we refer the reader to \cite{Cvet2004,Hum1972}; here, we only require the fact that both $E_6$ and $E_7$ are proper subsets of $E_8$ (and so $l(E_6)=l(E_7)=l(E_8)=\{\sqrt{2}\}$).

Let $\Sigma$ and $\Sigma'$ be two sets of vectors in  two Euclidean spaces, respectively. Then,  $\Sigma$ and $\Sigma'$ are said to be \textit{isometric} if their Gram matrices coincide, that is, there exists a bijection $\phi\colon\Sigma\to \Sigma'$ such that
\[
\langle \phi(\alpha), \phi(\beta) \rangle=\langle \alpha, \beta \rangle, \quad \text{for } \alpha,\beta\in \Sigma.
\]

As noted in Remark \ref{remark1}, an isomorphism between two irreducible root systems $\Phi,\Phi'$ is not required to preserve inner products. However, we can show that it scales the inner product by a  real constant.

\begin{proposition}\label{lemlen}
    Suppose that two root systems $\Phi$ and $\Phi'$ are isomorphic via an isomorphism~$\phi$. If~$\Phi$ is irreducible, then the following statements hold:
    \begin{enumerate}
        \item[(i)] There exists a constant $r> 0$ such that \[\langle\phi(\alpha),\phi(\beta)\rangle=r\langle\alpha,\beta\rangle, \quad\text{for any } \alpha,\beta\in\Phi.\]
        \item[(ii)] $l(\Phi')=\sqrt{r} l(\Phi)$.
        \item[(iii)] If $r=1$, then $\Phi$ is isometric to $\Phi'$.
    \end{enumerate}
\end{proposition}
\begin{proof} 
Part~(i) is the standard fact that an isomorphism between irreducible root systems is a similarity of the underlying Euclidean spaces; equivalently,
the invariant inner product is unique up to a positive scalar on each irreducible component. It can be deduced from a more general result established in~\cite[Chapter~VI, Section~1, no.~2, Proposition~7]{Bou2002}.

It remains to prove the latter two parts. By part~(i), we obtain $$l(\Phi')=\{\|\phi(\alpha)\|\,:\,\alpha\in\Phi\}=\{\sqrt{r}\|\alpha\|\,:\,\alpha\in\Phi\}=\sqrt{r}\cdot l(\Phi),$$ and so (ii) holds. If $r=1$, then $\Phi$ is isometric to $\Phi'$ with respect to the bijection  $\phi$ restricted to $\Phi$. Hence, (iii) holds.
\end{proof}

The next result is used only once, namely in the proof of Lemma~\ref{l52}, where it restricts the possible irreducible root systems containing representations with two root lengths.

\begin{proposition}\label{lemmet}
    If $\Phi$ is an irreducible root system with $l(\Phi)=\{1,\sqrt{2}\}$, then $\Phi$ is isometric to either $B_n$, or $\frac{1}{\sqrt{2}}C_n$ $($for some $n\ge 2)$, or $F_4$.
\end{proposition}
\begin{proof}
By Theorem~\ref{rootthm}, the root system $\Phi$ is isomorphic to one of
the irreducible root systems listed there. Under an isomorphism of
irreducible root systems, Proposition~\ref{lemlen}(ii) shows that all
root lengths are multiplied by a common positive constant. Consequently, the ratio of the root lengths is preserved.

Since $l(\Phi)=\{1,\sqrt2\}$, the root system $\Phi$ cannot be
isomorphic to $A_n$, $D_n$, $E_6$, $E_7$, $E_8$, or $G_2$, because
\[
    l(A_n)=l(D_n)=l(E_6)=l(E_7)=l(E_8)=\{\sqrt2\},~ l(G_2)=\{\sqrt2,\sqrt6\}.
\]
Therefore, $\Phi$ is isomorphic to one of $B_n$, $C_n$ $(n\ge2)$, or
$F_4$.

The standard root systems $B_n$ and $F_4$ have length set
$\{1,\sqrt2\}$, while $C_n$ has length set $\{\sqrt2,2\}$. Since
scaling does not change the isomorphism type of a root system, we replace
$C_n$ by $\frac1{\sqrt2}C_n$, whose length set is also
$\{1,\sqrt2\}$. Applying Proposition~\ref{lemlen}(ii) to an
isomorphism from $\Phi$ to one of these normalized systems yields a
scaling factor $r=1$. Hence, Proposition~\ref{lemlen}(iii) implies that
$\Phi$ is isometric to $B_n$, $\frac1{\sqrt2}C_n$, or $F_4$.
\end{proof}

\section{Representations of marked graphs}\label{sec4}

Let \(\mathbf G=(G,\mathbf b)\) be a marked graph. Suppose that
\(\mathbf G\) has a representation
\(\Sigma=(\alpha_1,\alpha_2,\ldots,\alpha_n)\). By~the identity \eqref{41e},
\[
    \gram(\Sigma)=\hat{A}(\mathbf G)+2I_n
    =A(G)+\mathbf b\mathbf b^{\top}-J_n+2I_n .
\]
It follows that, for \(1\le i\le n\),
\begin{equation}\label{42e}
    \langle\alpha_i,\alpha_i\rangle\in\{1,2\},
\end{equation}
and \(\langle\alpha_i,\alpha_i\rangle=2\) precisely when
\(i\in V_{\mathbf b}(G)\). We set
\[
    \Sigma_l=\{\alpha_i\in\Sigma:\|\alpha_i\|=\sqrt2\},
    \quad
    \Sigma_s=\{\alpha_i\in\Sigma:\|\alpha_i\|=1\},
\]
so that \(\Sigma=\Sigma_l\cup\Sigma_s\)  (disjointly). We refer to the vectors in $\Sigma_l$ and $\Sigma_s$ as the
\emph{long vectors} and \emph{short vectors} of the representation,
respectively. By \eqref{42e}, long vectors correspond to marked vertices, while short vectors correspond to unmarked vertices.

For \(i\ne j\), the same identity gives
\begin{equation}\label{con:innerproduct}
  \langle \alpha_i,\alpha_j\rangle\in
\begin{cases}
\{0,1\}, & \alpha_i,\alpha_j\in\Sigma_l,\\
\{0,-1\},& \alpha_i\in\Sigma_s\text{ or }\alpha_j\in\Sigma_s .
\end{cases}
\end{equation}

Conversely, let $\Sigma=(\alpha_1,\alpha_2,\ldots,\alpha_n)$ be a tuple of real vectors
satisfying \eqref{42e} and \eqref{con:innerproduct}. Define a binary vector
$\mathbf b=(b_i)_{i=1}^n$ by
$b_i=\langle \alpha_i,\alpha_i\rangle-1$, and set
$A=\operatorname{Gram}(\Sigma)-\mathbf b\mathbf b^\top+J_n-2I.$
The matrix $A$ is a $(0,1)$-matrix with zero
diagonal, and therefore is the adjacency matrix of a graph $G$.
It follows that $\Sigma$ is a representation of the marked graph $\mathbf G=(G,\mathbf b)$. To emphasize the role of the representation, we will also write $\mathbf G(\Sigma)$ for the marked graph $\mathbf G$ arising from~$\Sigma$.

\begin{figure}[h!]
  \centering
  \includegraphics[page=4,scale=0.9]{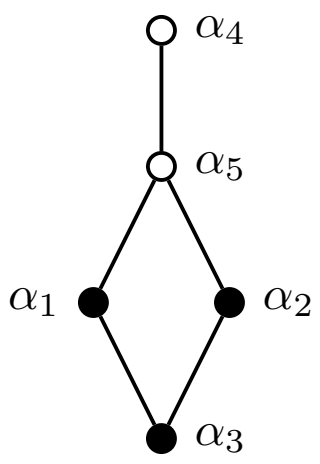}
  \caption{The  marked graph $\mathbf{G}(\Sigma)$ of Example~\ref{exa:e}.}
  \label{fig:5}
\end{figure}

\begin{example}\label{exa:e} Let
\begin{equation*}
	\Sigma=(\alpha_1,\alpha_2,\alpha_3,\alpha_4,\alpha_5)=\begin{pmatrix}
		1 & 1 & 1 & -1 & 0\\
		1 & -1 & 0 & 0 & 0\\
		0 & 0 & 1 & 0 & -1
	\end{pmatrix},
\end{equation*}
then $\mathbf{G}(\Sigma)$ is the  marked graph illustrated in Figure~\ref{fig:5}. Moreover, a direct algebraic computation confirms that $\mathbf{G}(\Sigma)$ is isomorphic to $\mathbf{L}^*(K_2(1,0))$.\end{example}

We now state several auxiliary lemmas.
\begin{lemma}\label{l41}
    A   marked graph $\mathbf{G}$ satisfies $ \lambda_{\min}(\hat{A}(\mathbf{G}))\ge -2$ if and only if $\mathbf{G}$ has a representation.
\end{lemma}
\begin{proof}
    If $\Sigma$ is a representation of $\mathbf{G}$, then $\gram(\Sigma)=\hat{A}(\mathbf{G})+2I$ is positive semidefinite, and so \( \lambda_{\min}\bigl(\hat{A}(\mathbf{G})\bigr) \geq -2\).

    If \( \lambda_{\min}\bigl(\hat{A}(\mathbf{G})\bigr) \geq -2\), then there exists a real  matrix $R$ such that $\hat{A}(\mathbf{G})+2I=R^\top R=\gram(\Sigma)$, where $\Sigma$ is the  set of   column vectors of $R$. Hence,  $\Sigma$ is  a representation of $\mathbf{G}$.
\end{proof}

The next lemma describes a deeper interplay between a  marked graph and its representation.

\begin{lemma}\label{l42}
Let \(\Sigma=(\alpha_1,\alpha_2,\ldots,\alpha_n)\) be a representation of a marked
graph \(\mathbf G=(G,\mathbf b)\) in a Euclidean space. Then:
\begin{enumerate}[label=(\roman*)]
    \item[(i)] \(\mathbf G\) is irreducible if and only if \(\Sigma\) is irreducible.
         \item[(ii)] If \(\sigma\) is an orthogonal transformation of the ambient 
    space, then \(\sigma(\Sigma)=(\sigma\alpha_1,\sigma\alpha_2,\ldots,\sigma\alpha_n)\) is also a
    representation of \(\mathbf G\).
     \item[(iii)] For each \(S\subseteq\{1,2,\ldots,n\}\), the marked graph
    \(\mathbf G(\Sigma_S)\), where \(\Sigma_S=(\alpha_i)_{i\in S}\), is an
    induced marked subgraph of \(\mathbf G=\mathbf G(\Sigma)\).
\end{enumerate}
\end{lemma}

\begin{proof}
(i): A partition \(S\cup S^c=\{1,2,\ldots,n\}\), with both parts
non-empty, gives a non-trivial marked-union decomposition of \(\mathbf G\) if
and only if the off-diagonal block of \(\hat{A}(\mathbf G)\) indexed by
\(S\times S^c\) is zero. Since
\[
    \gram(\Sigma)=\hat{A}(\mathbf G)+2I,
\]
the same off-diagonal block is zero if and only if
\(\langle \alpha_i,\alpha_j\rangle=0\) for all \(i\in S\) and
\(j\in S^c\). This is precisely a non-trivial orthogonal decomposition of
\(\Sigma\). Hence \(\mathbf G\) is irreducible if and only if \(\Sigma\) is
irreducible.

(ii): Since \(\sigma\) preserves inner products,
\(\gram(\sigma(\Sigma))=\gram(\Sigma)=\hat{A}(\mathbf G)+2I\). Hence \(\sigma(\Sigma)\) is
a representation of \(\mathbf G\).

(iii): By the construction of \(\mathbf G(\Sigma)\), passing from
\(\Sigma\) to \(\Sigma_S\) amounts to taking the principal submatrix indexed
by \(S\). Hence, \(\mathbf G(\Sigma_S)=(G[S],\mathbf b|_S)\), where \(G[S]\)
is the subgraph of \(G\) induced by \(S\) and $\mathbf{b}|_S$ is the restriction of $\mathbf{b}$ to $S$. Therefore, \(\mathbf G(\Sigma_S)\)
is an induced marked subgraph of \(\mathbf G\).
\end{proof}

Remark~\ref{rem:repT2} shows that, for a \(B\)-graph \(H\), the marked
graph \(\mathbf L^*(H)\) is represented by \(B_n\) for some
\(n\ge2\). We now prove that, conversely, every representation of a marked
graph is contained in a root system. We shall use the following standard
notion.

A subset \( S\) of a Euclidean space \(E\) is called
a \emph{lattice} if:
(a)  \(\alpha\pm\beta\in S\) whenever \(\alpha,\beta\in S\);
(b) \( S\) spans \(E\);
(c)  there exists \(r>0\) such that \(0\) is the only vector of
    \( S\) in the open ball \(\{x\in E:\|x\|<r\}\).

\begin{lemma}[\cite{Cas1997}]\label{lem:finite}
Let \( S\) be a lattice in a Euclidean space \(E\).
Then, for every \(r>0\), the set $\{\alpha\in S:\|\alpha\|<r\}$ is finite.
\end{lemma}

The following theorem forms the bridge between marked graph representations and  root system theory and is the key tool in the proof of the main results.

\begin{theorem}\label{lem:contained_in_rootsystem}
Let \(\Sigma=(\alpha_1,\alpha_2,\ldots,\alpha_n)\) be a representation of a marked
graph \(\mathbf G\) in a Euclidean space. Then \(\Sigma\) is contained in
a root system. Moreover, if \(\Sigma\) is irreducible, then \(\Sigma\) is
contained in an irreducible root system.
\end{theorem}

\begin{proof}
Let \(E=\operatorname{span}_{\mathbb R}\Sigma\) and 
\[
     S=\sum_{i=1}^n\mathbb Z\alpha_i
    =\left\{\sum_{i=1}^n z_i\alpha_i:z_i\in\mathbb Z\right\}\subseteq E .
\]
For \(x=\sum_i a_i\alpha_i\) and \(y=\sum_j b_j\alpha_j\) in \( S\),
where \(a_i,b_j\in\mathbb Z\), we have
$
    \langle x,y\rangle
    =\sum_{i,j}a_i b_j\langle\alpha_i,\alpha_j\rangle .
$
By \eqref{42e} and \eqref{con:innerproduct}, all inner products
\(\langle\alpha_i,\alpha_j\rangle\) are integers. Hence,
\begin{equation}\label{cla:integrality}
    \langle x,y\rangle\in\mathbb Z,\quad x,y\in S .
\end{equation}

We first show that \( S\) is a lattice in \(E\). It is closed under
addition and subtraction by definition, and it spans \(E\) since
\(\Sigma\subseteq S\). If \(x\in S\) and \(\|x\|<1\), then
\(\langle x,x\rangle\) is a non-negative integer smaller than \(1\), by
\eqref{cla:integrality}; hence \(x=0\). Thus \( S\) is a lattice.

Let $\Phi=\{\alpha\in  S\,:\,\langle \alpha,\alpha\rangle\in\{1, 2\}\}$. We  show that $\Phi$ is a root system in $E$ by verifying the defining rules (R1)--(R4).

(R1): Since $\Sigma\subseteq\Phi\subseteq  S$ and $\Sigma, S$ both span $E$, we obtain $\sspan_{\mathbb{R}}\Phi = E$. Observing that $\Phi$ is a subset of $ S\cap \{\alpha\in E\,:\,\|\alpha\|<\sqrt{2}+\varepsilon\}$ for some $\varepsilon>0$, and employing Lemma~\ref{lem:finite}, we deduce that $\Phi$ is finite.

(R2): Consider  $\alpha \in \Phi$. Since $-\alpha\in  S$ and $\langle-\alpha,-\alpha\rangle\in\{1,2\}$, we obtain $-\alpha\in\Phi$. Suppose that $c\alpha\in\Phi$ for some $c\ne 0$. Since $\langle\alpha,\alpha\rangle\in\{1,2\}$ and $c^2\langle\alpha,\alpha\rangle=\langle c\alpha,c\alpha\rangle\in\{1,2\}$, we obtain $c\in\{\pm\frac{1}{\sqrt{2}},\pm1,\pm\sqrt{2}\}$. Combining this with \eqref{cla:integrality}, we obtain $c\langle\alpha,\alpha\rangle=\langle c\alpha,\alpha\rangle\in\mathbb{Z}$ and so $c\in\{\pm1\}$.
Therefore, (R2) holds.

(R4): Consider  $\alpha, \beta \in \Phi$. Since $\langle\alpha, \alpha\rangle \in \{1, 2\}$ and $\langle\alpha, \beta\rangle \in \mathbb{Z}$, by \eqref{cla:integrality} we have
$
    \frac{2\langle\beta, \alpha\rangle}{\langle\alpha, \alpha\rangle} \in \mathbb{Z},
$
and so (R4) holds.

(R3): As before, consider $\alpha,\beta \in \Phi$. We need to show that
\[
\sigma_\alpha(\beta) = \beta - \frac{2\langle\beta, \alpha\rangle}{\langle\alpha, \alpha\rangle} \alpha\in\Phi.
\]
A direct algebraic computation leads to
\[
\langle \sigma_\alpha(\beta),\sigma_\alpha(\beta)\rangle = \Big\langle\beta - \frac{2\langle\beta, \alpha\rangle}{\langle\alpha, \alpha\rangle} \alpha,\beta - \frac{2\langle\beta, \alpha\rangle}{\langle\alpha, \alpha\rangle} \alpha\Big\rangle=\langle\beta,\beta\rangle\in\{1,2\}.
\]
Since $\frac{2\langle\beta, \alpha\rangle}{\langle\alpha, \alpha\rangle}\in\mathbb{Z}$ by  (R4), we have
\[
    \sigma_\alpha(\beta)=\beta - \frac{2\langle\beta, \alpha\rangle}{\langle\alpha, \alpha\rangle} \alpha\in  S.
\]
Hence, $\sigma_\alpha(\beta)\in\Phi$.

Therefore, \(\Phi\) is a  root system in \(E\). Suppose that $\Sigma$ is irreducible. By Lemma \ref{lemire}, $\Phi$ has a decomposition $\Phi=\bigcup_{i=1}^{r}\Phi_i$, where $r\ge 1$, $E=\bigoplus_{i=1}^{r}E_{i}$ is an orthogonal direct sum,  and each $\Phi_i$ is an irreducible root system in the subspace $E_{i}$. Note that for  $i\ne j$, every vector in $\Phi_i$ is orthogonal to every vector in $\Phi_j$. If \(\Sigma\) met two distinct components of this decomposition, then the corresponding non-empty parts of \(\Sigma\) would be mutually orthogonal, giving a non-trivial orthogonal decomposition of~\(\Sigma\). This contradicts the irreducibility of \(\Sigma\). Hence, \(\Sigma\subseteq \Phi_i\) for some \(1\le i\le r\), and so \(\Sigma\) is contained in an irreducible root system.
\end{proof}

\section{Proof of main result}\label{sec5}

We begin with the degenerate case when no vertices are marked. In this situation, the condition $\lambda_{\min}(\widehat A(G))\ge -2$ admits a  simple characterization.

\begin{lemma}\label{l51}
    Let $\mathbf{G}=(G,\mathbf{0})$ be an irreducible  marked graph. Then the following statements are equivalent:
    \begin{enumerate}
        \item[(i)] $ \lambda_{\min}(\hat{A}(\mathbf{G}))\ge -2$;
        \item[(ii)] $G$ is either $K_1$ or $2K_1$;
        \item[(iii)] $\mathbf{G}$ is represented by $B_2$.
    \end{enumerate}
\end{lemma}
\begin{proof}
    We establish the equivalence by proving the following cyclic implications.

    \noindent (i) $\Rightarrow$ (ii): Since $\mathbf{b}=\mathbf{0}$, the marked adjacency matrix satisfies
    \[
        \hat{A}(\mathbf{G})=A(G)-J=-A(\overline{G})-I .
    \]
    Hence,
    \[
         \lambda_{\min}(\hat{A}(\mathbf{G}))=-\lambda_1(A(\overline{G}))-1,
    \]
    and so the condition $ \lambda_{\min}(\hat{A}(\mathbf{G})) \ge -2$ is equivalent to $ \lambda_1(A(\overline{G})) \le 1$. By the classical spectral characterization of graphs with largest eigenvalue at most $1$ (see Section~3.11 in \cite{Cvet2010}), $\overline{G}$ is a disjoint union of copies of $K_1$ and $K_2$, or equivalently an induced subgraph of $tK_2$ for some $t \ge 1$. Taking complements, $G$ is therefore an induced subgraph of the cocktail party graph $CP(t)=\overline{tK_2}$.
    
    Such an induced subgraph is a complete multipartite graph whose parts have size at most two. If at least two parts are non-empty, then $G$ is a non-trivial join of two induced subgraphs. Since all vertices are unmarked, this gives a non-trivial marked-union decomposition of $\mathbf{G}$, contradicting irreducibility. Thus only one part is non-empty, and hence $G$ is either $K_1$ or $2K_1$.

    \noindent (ii) $\Rightarrow$ (iii): If $G=K_1$, then the \(1\)-tuple $(\mathbf{e}_1)$ is a representation of $\mathbf{G}$ in $B_2$. If $G=2K_1$, then the \(2\)-tuple $(\mathbf{e}_1,-\mathbf{e}_1)$ is a representation of $(2K_1,\mathbf{0})$ in $B_2$.

    \noindent (iii) $\Rightarrow$ (i): This implication follows immediately from Lemma~\ref{l41}.
\end{proof}

We now turn to irreducible marked graphs with at least one marked vertex. The next lemma shows that every such graph satisfying $\lambda_{\min}(\widehat A(G))\ge -2$ is represented by either $B_n$ or $E_8$.

\begin{lemma}\label{l52}
    Let $\mathbf{G}=(G,\mathbf{b})$ be an irreducible  marked graph with $\mathbf{b}\ne \mathbf{0}$. Then $ \lambda_{\min}(\hat{A}(\mathbf{G}))\ge -2$ if and only if $\mathbf{G}$ is represented by either $E_8$ or some $B_n$.
\end{lemma}
\begin{proof}
The sufficiency follows immediately from Lemma \ref{l41}.

Next, we prove the necessity. The argument proceeds by considering separately the cases
$\mathbf b=\mathbf j$ and $\mathbf b\neq \mathbf j$.
In the latter case, the representation theory developed in the previous section, together with Proposition~\ref{lemmet}, reduces the problem to the root systems $B_n$, $\frac1{\sqrt2}C_n$, and $F_4$; we then show that the latter two possibilities always yield a representation in $B_n$.

We first consider the case  $\mathbf{b}=\mathbf{j}$.
In this case \(\hat{A}(\mathbf G)=A(G)\), and the present assumption gives
\( \lambda_{\min}(A(G))\ge -2\). Moreover, by the observation following the
definition of irreducibility, \(G\) is connected. By Theorem \ref{thm:Cam}
and \(D_r\subset B_r\) for all admissible \(r\), we obtain that
\(\mathbf{G}\) is represented by either \(E_8\) or some \(B_n\).

We next consider the case $\mathbf{b}\ne\mathbf{j}$ (and $\mathbf{b}\ne\mathbf{0}$). 
By Lemmas \ref{l41} and \ref{l42}(i), we obtain that $\mathbf{G}$ has a  representation, and every representation of  $\mathbf{G}$ is irreducible. Let $\Sigma$ be such a representation. By Theorem~\ref{lem:contained_in_rootsystem}, $\Sigma$ is contained in an irreducible root system $\Phi$, and $\Sigma$ satisfies conditions (\ref{42e}) and (\ref{con:innerproduct}). 

Since $\mathbf{b}\ne \mathbf{0}$ and $\mathbf{b}\ne \mathbf{j}$, we obtain $\Sigma_l\ne \emptyset$ and $\Sigma_s\ne \emptyset$, and further \( l(\Phi) = \{1, \sqrt{2}\} \) by Lemma \ref{lem:rootlength}. By Proposition~\ref{lemmet}, $\Phi$ is isometric to either $B_n$, or
$\frac{1}{\sqrt{2}}C_n$ for some $n\ge 2$, or $F_4$. Replacing $\Sigma$ by its image under this isometry, which preserves the Gram matrix, we may assume that $\Sigma$ is contained in one of these standard root systems, and  therefore it suffices to consider the cases $\Sigma\subseteq \frac{1}{\sqrt{2}}C_n$ and $\Sigma\subseteq F_4$, and prove that in both situations $G$ admits a representation in $B_n$.

\smallskip\noindent\textit{Case 1: $\Sigma\subseteq \frac{1}{\sqrt{2}}C_n$.} 
Here we have
\[
\emptyset\ne\Sigma_l\subseteq\left\{\pm\sqrt{2}\mathbf{e}_i\,:\,1\le i\le n\right\}
\quad\text{and}\quad
\emptyset\ne\Sigma_s\subseteq\left\{\tfrac{1}{\sqrt{2}}(\pm\mathbf{e}_i\pm\mathbf{e}_j)\,:\,1\le i<j\le n\right\}.
\]
Since $\Sigma_s\ne \emptyset$, pick
\[
v_1=\frac{1}{\sqrt{2}}(\varepsilon_p\mathbf{e}_p+\varepsilon_q\mathbf{e}_q)\in\Sigma_s,
\]
where $p\ne q$ and $\varepsilon_p,\varepsilon_q\in\{\pm1\}$. Let $\tau$ be a signed permutation of the standard basis such that
\[
\tau(\varepsilon_p\mathbf{e}_p)=\mathbf{e}_1\quad \text{and}\quad
\tau(\varepsilon_q\mathbf{e}_q)=\mathbf{e}_2.
\]
Then $\tau$ is an orthogonal transformation preserving $\frac{1}{\sqrt{2}}C_n$. By Lemma~\ref{l42}(ii), \(\tau(\Sigma)\) is also a representation of \(\mathbf G\). Replacing \(\Sigma\) by \(\tau(\Sigma)\) and \(v_1\) by \(\tau(v_1)\), we may assume
$
v_1=\frac{1}{\sqrt{2}}(\mathbf{e}_1+\mathbf{e}_2).
$
By \eqref{con:innerproduct}, we obtain
\[
\frac{1}{\sqrt{2}}(\pm\mathbf{e}_1\pm\mathbf{e}_k),\,
\frac{1}{\sqrt{2}}(\pm\mathbf{e}_2\pm\mathbf{e}_k)\not\in \Sigma,
\quad\text{for } k\not\in\{1,2\}.
\]
Thus $\emptyset\ne\Sigma_s\subseteq \{\tfrac{1}{\sqrt{2}}(\pm\mathbf{e}_1\pm\mathbf{e}_2)\}.$ Since $\Sigma$ is irreducible, it follows that
\[
\Sigma_l\subseteq\{\pm\sqrt{2}\mathbf{e}_1,\pm\sqrt{2}\mathbf{e}_2\}.
\]

Let $v_2=\frac{1}{\sqrt{2}}(\mathbf{e}_1-\mathbf{e}_2)$. Note that $\sqrt{2}\mathbf{e}_1=v_1+v_2$, $\sqrt{2}\mathbf{e}_2=v_1-v_2$, and $\{v_1,v_2\}$ forms an orthonormal basis of $\mathbb{R}^2$. Let $\sigma$ be an orthogonal transformation on $\mathbb{R}^2$ such that $\sigma(v_k)=\mathbf{e}_k$ for $k\in\{1,2\}$. Since
\[
\Sigma\subseteq\{\pm\sqrt{2}\mathbf{e}_1,\pm\sqrt{2}\mathbf{e}_2,\tfrac{1}{\sqrt{2}}(\pm\mathbf{e}_1\pm\mathbf{e}_2)\}
=\{\pm v_1\pm v_2,\pm v_1,\pm v_2\},
\]
we obtain
$\sigma(\Sigma)\subseteq\{\pm \mathbf{e}_1\pm \mathbf{e}_2,\pm \mathbf{e}_1,\pm \mathbf{e}_2\}=B_2$.
By Lemma \ref{l42}(ii), $\mathbf{G}$ is represented by $B_2$.

\smallskip\noindent\textit{Case 2: $\Sigma\subseteq F_4$.} 
Now, we have \[\emptyset\ne\Sigma_l\subseteq \{\pm\mathbf{e}_i\pm\mathbf{e}_j\,:\,1\le i<j\le 4\}\text{ and }\emptyset\ne\Sigma_s\subseteq \{\pm\mathbf{e}_i\,:\,1\le i\le 4\}\cup \left\{\tfrac{1}{2}(\pm\mathbf{e}_1\pm\mathbf{e}_2\pm\mathbf{e}_3\pm\mathbf{e}_4)\right\}.\]
If $\Sigma_s\cap \left\{\tfrac{1}{2}(\pm\mathbf{e}_1\pm\mathbf{e}_2\pm\mathbf{e}_3\pm\mathbf{e}_4)\right\}=\emptyset $, then $\Sigma\subseteq B_4$, and there is nothing to prove.

Henceforth, we  suppose that $\Sigma_s\cap \left\{\tfrac{1}{2}(\pm\mathbf{e}_1\pm\mathbf{e}_2\pm\mathbf{e}_3\pm\mathbf{e}_4)\right\}\ne \emptyset$. Note that $\langle\alpha,\beta\rangle\in\{\pm 1/2\}$ for  $\alpha\in \{\pm\mathbf{e}_i\,:\,1\le i\le 4\}$ and $\beta\in \left\{\tfrac{1}{2}(\pm\mathbf{e}_1\pm\mathbf{e}_2\pm\mathbf{e}_3\pm\mathbf{e}_4)\right\}$. By \eqref{con:innerproduct}, we obtain
$\emptyset\ne\Sigma_s\subseteq \left\{\tfrac{1}{2}(\pm\mathbf{e}_1\pm\mathbf{e}_2\pm\mathbf{e}_3\pm\mathbf{e}_4)\right\}$.

Since \(\Sigma_s \neq \emptyset\), pick
\[
v_1=\frac{1}{2}(\varepsilon_1\mathbf{e}_1+\varepsilon_2\mathbf{e}_2+\varepsilon_3\mathbf{e}_3+\varepsilon_4\mathbf{e}_4)\in\Sigma_s,
\]
where \(\varepsilon_i\in\{\pm1\}\). Let \(\sigma\) be the diagonal signed orthogonal transformation defined by \(\sigma(\varepsilon_i\mathbf e_i)=\mathbf e_i\) for \(1\le i\le 4\). Then \(\sigma\) preserves \(F_4\). By Lemma~\ref{l42}(ii), \(\sigma(\Sigma)\) is also a representation of \(\mathbf G\). Replacing \(\Sigma\) by  \(\sigma(\Sigma)\)  and \(v_1\) by \(\sigma(v_1)\), we may assume
\[
v_1 = \frac{1}{2}(\mathbf{e}_1+\mathbf{e}_2+\mathbf{e}_3+\mathbf{e}_4).
\]
Let \(\omega = \frac{1}{2}(\delta_1\mathbf{e}_1+\delta_2\mathbf{e}_2+\delta_3\mathbf{e}_3+\delta_4\mathbf{e}_4) \in \Sigma_s\setminus\{\pm v_1\}\), where \(\delta_i \in\{ \pm 1\}\). By \eqref{con:innerproduct}, we obtain \(\langle v_1, \omega \rangle \in\{0,-1\}\). Since both \(v_1\) and \(\omega\) are unit vectors and \(\omega\ne -v_1\), the value \(-1\) is impossible. Therefore, \(\langle v_1, \omega \rangle = \frac{1}{4}(\delta_1 + \delta_2 + \delta_3 + \delta_4) = 0\). Hence, \(\{\delta_i\}\) has exactly two \(+1\)s and two \(-1\)s. Up to sign, the only such vectors are
\[
v_2 = \frac{1}{2}(-\mathbf{e}_1+\mathbf{e}_2+\mathbf{e}_3-\mathbf{e}_4),\quad
v_3 = \frac{1}{2}(-\mathbf{e}_1+\mathbf{e}_2-\mathbf{e}_3+\mathbf{e}_4),\quad
v_4 = \frac{1}{2}(-\mathbf{e}_1-\mathbf{e}_2+\mathbf{e}_3+\mathbf{e}_4).
\]
Thus, \(\Sigma_s \subseteq \{\pm v_1, \pm v_2, \pm v_3, \pm v_4\}\).

Note that \(\{v_k\}_{k=1}^4\) forms an orthonormal basis of \(\mathbb{R}^4\). Let \(\tau\) be an orthogonal transformation on~$\mathbb{R}^4$ such that \(\tau(v_k) = \mathbf{e}_k\) for \(k=1,2,3,4\). Then \(\tau(\Sigma_s) \subseteq \{\pm \mathbf{e}_k : 1\le k\le 4\} \subset B_4\).

Now let \(\alpha \in \Sigma_l\). Then \(\alpha\) is of the form \(\pm \mathbf{e}_i \pm \mathbf{e}_j\) for some \(1 \le i < j \le 4\). This implies that \(c_k:=\langle \alpha, v_k \rangle\in \{0, \pm 1\}\) for $ k\in\{1,2,3,4\}$.  Since \(\{v_k\}_{k=1}^4\) forms an orthonormal basis of \(\mathbb{R}^4\), we have
\[
\alpha = \sum_{k=1}^4 c_k v_k, \quad \text{and}\quad \langle \alpha,\alpha\rangle=\sum_{k=1}^4 c_k^2 = 2.
\]
Hence, exactly two of the \(c_k\) are non-zero, say \(c_r\) and \(c_s\) with \(r \neq s\) and \(c_r, c_s \in \{\pm 1\}\). It follows that $\alpha = c_r v_r + c_s v_s,$
and therefore
\[
\tau(\alpha) = c_r \mathbf{e}_r + c_s \mathbf{e}_s \in \{\pm \mathbf{e}_r \pm \mathbf{e}_s\} \subset B_4.
\]

Consequently, \(\tau(\Sigma) \subset B_4\). By Lemma \ref{l42}(ii), $\mathbf{G}$ is represented by $B_4$.
\end{proof}

The next step is to characterize those irreducible marked graphs that admit a representation in~$B_n$. This leads to the marked analogue of generalized line graphs.

\begin{lemma}\label{l53}
    Let $\mathbf{G}=(G,\mathbf{b})$ be an irreducible marked graph with $\mathbf{b}\ne \mathbf{0}$. Then $\mathbf{G}$ is represented by some $B_n$ if and only if $\mathbf{G}$ is an induced subgraph of an extended generalized line graph.
\end{lemma}
\begin{proof} Suppose that $\mathbf{G}$ is an induced subgraph of an extended generalized line graph. By Remark~\ref{rem:repT2}, every extended generalized line graph is represented by some $B_n$. Restricting the corresponding representation to the vertices of the induced subgraph, we obtain that $\mathbf{G}$ is represented by~$B_n$.

Next, we prove the necessity. The proof proceeds by analyzing the structure of a representation $\Sigma\subseteq B_n$. We first partition the coordinates according to the signs occurring in the long vectors of $\Sigma$, and then show that every vector of $\Sigma$ belongs to a distinguished set of roots naturally associated with a suitable $B$-graph. This will imply that $G$ is an induced subgraph of the corresponding extended generalized line graph.

Suppose that $\Sigma\subseteq B_n$ is a representation of $\mathbf{G}$ for some $n\ge 2$. 
    Since $\mathbf{b}\ne \mathbf{0}$, we have
\begin{equation*}
    \emptyset\ne\Sigma_l\subseteq \{\pm\mathbf{e}_i\pm\mathbf{e}_j\,:\,1\le i<j\le n\} \text{ and }  \Sigma_s \subseteq \{\pm\mathbf{e}_i\,:\,1\le i\le n\}.
\end{equation*}
Let $N=\{1,2,\ldots,n\}$. We define $\supp(\mathbf{x})\coloneqq  \{ i \in N \,:\, x_i \neq 0 \}$ for  $\mathbf{x}=(x_1, x_2, \ldots, x_n)^\top$. After discarding unused coordinates if necessary, we may assume that 
$$\bigcup_{\alpha \in \Sigma} \supp(\alpha) = N.$$
We claim that
\begin{equation}\label{eq:long-support-cover}
    \bigcup_{\alpha \in \Sigma_l} \supp(\alpha) = N.
\end{equation}
To see this, suppose for contradiction that $\bigcup_{\alpha \in \Sigma_l} \supp(\alpha) \subsetneq N$. 
Then some coordinate \(p\in N\) appears only in vectors of \(\Sigma_s\), and hence  the vectors of \(\Sigma\) whose support contains $p$ would be orthogonal to the vectors of $\Sigma$ whose support does not contain $p$. Since \(\Sigma_l\ne\emptyset\), this would give a non-trivial orthogonal decomposition of \(\Sigma\), contradicting the irreducibility of \(\Sigma\). This contradiction shows that \eqref{eq:long-support-cover} holds.

For  $\alpha \in \Sigma$, we denote its $q$-th entry by $\alpha_q$. Then $\alpha_q\in\{0,\pm1\}$. We partition $N$ according to the signs occurring in each coordinate among the long vectors:
\begin{align*}
   & N_1\,\, \coloneqq  \{ q \in N \,:\, \exists\, \{\alpha, \beta\} \subset \Sigma_l \text{ such that } \alpha_q=1 \text{ and } \beta_q = -1 \},\\
   & N^+ \coloneqq  \{ q\in N \,:\, \forall \alpha\in \Sigma_l,\,\alpha_q\in\{0,1\}\},\\
   & N^- \coloneqq  \{ q\in N \,:\, \forall \alpha\in \Sigma_l,\,\alpha_q\in\{0,-1\}\}.
\end{align*}
If $N^- \neq \emptyset$, let $\sigma$ be the orthogonal transformation on $\mathbb{R}^n$ defined by
\[
 \sigma(\mathbf{e}_p) =
\begin{cases}
-\mathbf{e}_p, & \text{if } p \in N^-, \\
\mathbf{e}_p, & \text{otherwise}.
\end{cases}
\]
Since $\sigma(B_n) = B_n$, Lemma \ref{l42}(ii) implies that $\sigma(\Sigma) \subseteq B_n$ is also a representation of $\mathbf{G}$. Replacing $\Sigma$ by $\sigma(\Sigma)$ and redefining $N_1, N^+, N^-$ accordingly, we may assume $N^- = \emptyset$ for the rest of the proof.

Consider any $q\in N_1$. There exist  $\alpha,\beta\in \Sigma_l$ with $\alpha_q=1$ and $\beta_q=-1$. Write
\[
    \alpha=s_\alpha\mathbf e_p+\mathbf e_q,\quad
    \beta=s_\beta\mathbf e_k-\mathbf e_q,
\]
where \(s_\alpha,s_\beta\in\{\pm1\}\). Since \(\langle\alpha,\beta\rangle\in\{0,1\}\) by \eqref{con:innerproduct}, we must have \(p=k\) and \(s_\alpha=s_\beta\). Thus, for some \(s\in\{\pm1\}\),
\[
    \alpha=s\mathbf e_p+\mathbf e_q,\quad
    \beta=s\mathbf e_p-\mathbf e_q.
\]
Moreover, this pair is unique. Indeed, if \(\gamma\in\Sigma_l\) and \(\gamma_q=1\), then applying \eqref{con:innerproduct} to \(\gamma\) and \(\beta\) gives \(\gamma=\alpha\). Similarly, if \(\gamma_q=-1\), then \(\gamma=\beta\). Thus,
\begin{equation}\label{eqpetal}
	\text{$(\alpha,\beta)$ is the unique vector pair in $\Sigma_l$ such that $\alpha_q=1$ and $\beta_q=-1$.}
\end{equation}

We claim that \(p\notin N_1\). Otherwise, there exists \(\gamma\in\Sigma_l\) with \(\gamma_p=-s\). If \(\supp(\gamma)=\{p,r\}\) with \(r\ne q\), then \(\langle\gamma,\alpha\rangle=-1\), contradicting \eqref{con:innerproduct}. Hence, \(r=q\). But then \(\gamma_q=1\) or \(\gamma_q=-1\), contradicting the uniqueness in \eqref{eqpetal}, because \(\gamma_p=-s\) while \(\alpha_p=\beta_p=s\). Therefore, \(p\notin N_1\). Since \(p\) appears in a long vector and \(N^-=\emptyset\), we have \(p\in N^+\).
For each $q \in N_1$, let $u(q)$ denote the index $p$ determined by $q$ in the above argument. Then $u(q) \in N^+$ for every $q \in N_1$. Now, we obtain that
\[\Sigma_l\subseteq\{\mathbf{e}_{p_1}+\mathbf{e}_{p_2}\,:\,p_1\ne p_2\in N^+\}\cup\{\mathbf{e}_{u(q)}\pm\mathbf{e}_q\,:\,q\in N_1\}.\]

If $\alpha\in \Sigma_s$ with $\supp(\alpha)=\{p\}$, then \(p\notin N_1\). Indeed, if \(p\in N_1\), then there are long vectors \(s\mathbf e_r+\mathbf e_p\) and \(s\mathbf e_r-\mathbf e_p\) in \(\Sigma_l\); one of these has inner product \(1\) with \(\alpha=\mathbf e_p\) or with \(\alpha=-\mathbf e_p\), contradicting \eqref{con:innerproduct}. Hence, \(p\in N^+\). Since \(p\) appears in a long vector and \(N^-=\emptyset\), some long vector has \(p\)-th entry equal to \(1\). Therefore, \(\alpha\ne\mathbf e_p\), again by \eqref{con:innerproduct}, and so \(\alpha=-\mathbf e_p\). Thus,
\[\Sigma_s\subseteq\{-\mathbf{e}_p\,:\,p\in N^+\}.\]
Hence,
\[\Sigma\subseteq\Sigma^*\coloneqq \{\mathbf{e}_{p_1}+\mathbf{e}_{p_2}:p_1\ne p_2\in N^+\}\cup\{\mathbf{e}_{u(q)}\pm\mathbf{e}_q\,:\,q\in N_1\}\cup\{-\mathbf{e}_p\,:\,p\in N^+\}.\]

Let $t=|N^+|$. Since \(N\ne\emptyset\), \(N^-=\emptyset\), and every
coordinate in \(N_1\) has its associated coordinate \(u(q)\in N^+\), we have
\(t\ge1\). For each $p\in N^+$, let
$$a_p=|\{q\in N_1 \,:\, u(q)=p\}|.$$ Identifying the vertices of \(K_t\)
with the coordinates in \(N^+\), and the petal vertices attached at \(p\)
with those \(q\in N_1\) satisfying \(u(q)=p\), Remark~\ref{rem:repT2} gives
\[
    \mathbf{G}(\Sigma^*) \cong \mathbf{L}^*(H),
\]
where \(H=K_t(\mathbf a)\) and \(\mathbf a=(a_p)_{p\in N^+}\). Since \(\Sigma\subseteq\Sigma^*\), Lemma~\ref{l42}(iii) implies that $\mathbf{G}$ is an induced subgraph of $\mathbf{L}^*(H)$.
\end{proof}

It remains to identify the irreducible marked graphs represented by $E_8$ but not by any $B_n$.

\begin{lemma}\label{l54}
    Let $\mathbf{G}=(G,\mathbf{b})$ be an irreducible marked graph with $\mathbf{b}\ne \mathbf{0}$. Then $\mathbf{G}$ is represented by $E_8$ and is not represented by any $B_n$ if and only if $\mathbf{b}=\mathbf{j}$ and $G$ is an exceptional graph.
\end{lemma}
\begin{proof}
    Suppose that $\mathbf{G}$ is represented by $E_8$ and is not represented by any $B_n$. Since \(l(E_8)=\{\sqrt{2}\}\), the condition~\eqref{42e} implies that every entry of \(\mathbf{b}\) is equal to \(1\), and hence \(\mathbf{b}=\mathbf{j}\). Thus, \(\hat{A}(\mathbf G)=A(G)\). Moreover, by the observation following the definition of irreducibility, \(G\) is connected.
    
    Since \(D_n\subset B_n\) for every admissible \(n\), the graph \(G\) is not represented by any \(D_n\). By Theorem~\ref{thm:Cam}(iii), \(G\) is an exceptional graph.
    
    Conversely, suppose that $\mathbf{b}=\mathbf{j}$ and $G$ is an exceptional graph. Then $\hat{A}(\mathbf{G})=A(G)$ is represented by $E_8$ and is not represented by any  $D_n$, by Theorem~\ref{thm:Cam}(iii). If \(\mathbf G\) is represented by some \(B_n\), then, since all vertices are marked, every representing vector would have length \(\sqrt2\) by \eqref{42e}. Hence, the representing vectors would lie in the long-root subsystem \(D_n\subset B_n\), and \(G\) would be represented by \(D_n\), a contradiction. 
\end{proof}

Combining the previous lemmas, we obtain the announced characterization of irreducible marked graphs.

\medskip
\par\noindent\textit{Proof of Theorem $\ref{thm:-2rep}$.} 
Statement (i) follows from Lemma~\ref{l51} when \(\mathbf b=\mathbf{0}\), and from Lemma~\ref{l52} when \(\mathbf b\ne\mathbf{0}\). 
For statement (ii), the case \(\mathbf b=\mathbf{0}\) is covered by Lemma~\ref{l51}, while the case \(\mathbf b\ne\mathbf{0}\) follows from Lemma~\ref{l53}. 
Finally, statement (iii) follows from Lemma~\ref{l54}; the case \(\mathbf b=\mathbf{0}\) cannot occur here by Lemma~\ref{l51}. 
This completes the proof.
\qed
\medskip

The characterization of arbitrary marked graphs now follows from the irreducible decomposition and the corresponding characterization of irreducible marked graphs established in Theorem~\ref{thm:-2rep}.

\medskip

\par\noindent\textit{Proof of Theorem $\ref{thm2}$.} Two implications are considered separately.

\smallskip\noindent\textit{Sufficiency:} 
Let $\mathbf{G} = \widehat{\cup}_{i=1}^{r} \mathbf{G}_i$ be the irreducible decomposition and consider $\mathbf{G}_i=(G_i,\mathbf{b}_i)$. By Theorem~\ref{thm:-2rep}(i)--(iii), each  $\mathbf{G}_i$ of type (T1), (T2), or (T3) is represented by either $B_n$ for some \(n\ge2\) or by $E_8$. In either case, Theorem~\ref{thm:-2rep}(i) ensures $ \lambda_{\min}(\hat{A}(\mathbf{G}_i)) \ge -2$. 
Applying the block diagonal form in \eqref{eqq} yields
\[ \lambda_{\min}(\hat{A}(\mathbf{G}))=\min\{ \lambda_{\min}(\hat{A}(\mathbf{G}_i)):1\leq i\leq r\}\ge -2.\]

\smallskip\noindent\textit{Necessity:}
Assume that $\mathbf{G}$ is a   marked graph with $ \lambda_{\min}\bigl(\hat{A}(\mathbf{G})\bigr) \geq -2$. Then, $\mathbf{G}$ has the irreducible decomposition $\mathbf{G} = \widehat{\cup}_{i=1}^{r} \mathbf{G}_i$.
By the block diagonal form in \eqref{eqq},
\[
 \lambda_{\min}(\hat{A}(\mathbf{G}))=\min\{ \lambda_{\min}(\hat{A}(\mathbf{G}_i)):1\leq i\leq r\}.
\]
Thus, $ \lambda_{\min}(\hat{A}(\mathbf{G}_i))\ge -2$ for all $i$. By using Theorem~\ref{thm:-2rep}(i)--(iii) again, 
each $\mathbf{G}_i$ is necessarily of type (T1), (T2), or~(T3), and the proof is completed.
\qed

\section{Graphs whose second largest eigenvalue does not exceed 1}\label{sec6}

We conclude with an application of Theorem~\ref{thm2} to the study of graphs whose second largest eigenvalue does not exceed 1.

The problem of characterizing graphs whose second largest eigenvalue is at most 1 was posed by A.J.~Hoffman \cite{Cvet1982} in the early  1980s. Over the past  decades, this problem has attracted considerable attention from researchers and remains open in general. To date, it has been resolved only for certain special classes of graphs, including the classical cases of bipartite graphs~\cite{Pet1991}, line graphs~\cite{Pet1998}, and generalized line graphs~\cite{Pet2000}, as well as more recent results for outerplanar graphs~\cite{Li2023}, diamond-free graphs~\cite{Liu2024}, and claw-free graphs~\cite{Chen2025}. Further progress and related results can be found in \cite{Cheng2018, Cheng2025, Geng2024, Liu2021, Stanic2007, Stanic2009, Stanic2010}.

Let $\overline{G}$ denote the complement of a graph $G$. The study of graphs
with second largest eigenvalue at most $1$ is closely related to the
classical theory of graphs with least eigenvalue at least $-2$. 

\begin{theorem}[\cite{Cvet1982}]\label{thm1}
If a graph $G$ of order $n\ge 2$ satisfies
$ \lambda_{\min}(A(\overline{G})) \geq -2$, then $ \lambda_2(G)\leq 1$.
\end{theorem}

The following preliminary result is part of the folklore in spectral graph theory.

\begin{lemma}[{Weyl's inequalities, see \cite[Theorem~2.8.1]{Weyl1912}}]\label{weyl} Let \( M \) and \( N \) be \( n \times n \) Hermitian matrices, and let their eigenvalues be arranged in non-increasing order. Then,  for  \( 1\le i\le j\le n \),
\[
  \lambda_{i}(M+N)\ge  \lambda_{j}(M) +  \lambda_{n+i-j}(N).
\]
\end{lemma}

An $n$-vertex graph $G$ is said to be \textit{$-2$-markable} if there exists a vector $\mathbf{b} \in \{0,1\}^n$ such that $$ \lambda_{\min}(A(\overline{G})+\mathbf{b}\mathbf{b}^\top-J_n)\ge -2.$$ 
As a generalization of Theorem~\ref{thm1}, we show that every $-2$-markable graph $G$ satisfies $ \lambda_2(G) \leq 1$.

\begin{proposition}\label{thm3}
    If a graph $G$ of order $n\ge 2$ is $-2$-markable, then
$\lambda_2(G)\le 1$.
\end{proposition}
\begin{proof}
    Since $G$ is $-2$-markable, we may suppose that $ \lambda_{\min}(\hat{A}(\overline{G},\mathbf{b}))\ge -2$ for some binary vector $\mathbf{b}$. Note that
    \[
        \hat{A}(\overline{G},\mathbf{b})
        =A(\overline{G})+\mathbf{b} \mathbf{b}^\top-J_n
        =\mathbf{b} \mathbf{b}^\top-A(G)-I.
    \]
    Since \(\mathbf b\mathbf b^\top\) has rank at most one, we have
    $\lambda_2(\mathbf b\mathbf b^\top-I)=-1$.
    By inserting $i=2$, $j=n$, $M=\mathbf{b} \mathbf{b}^\top - A(G) - I$ and $N=A(G)$ in Lemma \ref{weyl}, we obtain
    \[-1= \lambda_2(\mathbf{b}\mathbf{b}^\top-I)\ge  \lambda_{\min}(\mathbf{b} \mathbf{b}^\top - A(G) - I)+ \lambda_2(G)\ge -2+ \lambda_2(G),\]
    which implies that $\lambda_2(G)\le 1$.
\end{proof}

The converse statement does not hold, but we provide some statistical data concerning graphs with 8 vertices.  There are exactly $2182$ connected graphs of order $8$ satisfying $\lambda_2 \leq 1$. Among them, $1731$ (close to 79.3\%) have complements with least eigenvalue at least $-2$, and $2132$ (close to 97.7\%) are $-2$-markable graphs; the details are given at the location specified in the Data availability section. These empirical observations indicate a  strong correlation between $-2$-markability and the spectral property $\lambda_2 \leq 1$.

By Theorem~\ref{thm2}, the marked graphs satisfying
$\lambda_{\min}\bigl(\hat{A}(\mathbf G)\bigr)\ge -2$
have been completely characterized in terms of the three types
(T1)--(T3). Therefore, a graph $G$ is $-2$-markable if and only if, for some binary vector $\mathbf b$, the
marked graph $(\overline G,\mathbf b)$ can be expressed as a marked union of some
irreducible marked graphs of types (T1)--(T3).

We next record a construction consequence in terms of marked joins.
For a marked graph $\mathbf{G}=(G,\mathbf{b})$, we denote its complement by $\overline{\mathbf{G}}=(\overline{G},\mathbf{b})$.
Let $\mathbf G_i=(G_i,\mathbf b_i)$, $1\leq i\leq r$, be pairwise
vertex-disjoint marked graphs. Note that the marked join and marked union are
complement-dual operations. In particular, for two marked graphs
$\overline{\mathbf G_1 \widehat{\vee} \mathbf G_2}=
 \overline{\mathbf G_1}\widehat{\cup}\overline{\mathbf G_2},$
and, more generally,
\begin{equation}\label{eq:vee}\overline{\widehat{\vee}_{i=1}^{r} \mathbf G_i}=
\widehat{\cup}_{i=1}^{r} \overline{\mathbf G_i}.\end{equation}

By \eqref{eq:vee}, saying that a marked graph $\mathbf G$ is obtained as a
marked join of $\mathbf G_1,\mathbf G_2,\ldots,\mathbf G_r$ is equivalent to saying that
$\overline{\mathbf G}$ is obtained as the marked union of
$\overline{\mathbf G}_1,\overline{\mathbf G}_2,\ldots,\overline{\mathbf G}_r$. Hence, the following corollary gives a useful
sufficient condition for producing graphs with second largest eigenvalue at
most $1$.

\begin{corollary}\label{cor2}
    Let $\{\mathbf{G}_i\}_{i=1}^r$ be a finite family  of pairwise vertex-disjoint marked graphs, and let $G$ be the underlying graph of $\mathbf{G}=\widehat{\vee}_{i=1}^{r} \mathbf{G}_i$. If, for $1\le i\le r$, $\overline{\mathbf{G}}_i$ is one of the types {\rm(T1)--(T3)} described in Theorem~{\rm \ref{thm2}}, then $\lambda_2(G)\le 1$.
\end{corollary}
\begin{proof} By Theorem~\ref{thm2} and \eqref{eq:vee}, we deduce
\[\lambda_{\min}(\hat{A}(\overline{\mathbf{G}}))
    =\lambda_{\min}\!\left(\hat{A}\!\left(\overline{\widehat{\vee}_{i=1}^{r}\mathbf{G}_i}\right)\right)
    =\lambda_{\min}\!\left(\hat{A}\!\left(\widehat{\cup}_{i=1}^{r}\overline{\mathbf{G}}_i\right)\right)\ge -2.\]
Thus, $G$ is $-2$-markable and the desired result follows from Proposition~\ref{thm3}.
\end{proof}

Corollary~\ref{cor2} can also be regarded as a generalization of \cite[Theorem~3.1]{Liu2025}. In what follows, we present two examples illustrating the use of the marked join operation to construct families of graphs whose second largest eigenvalue does not exceed 
$1$. 

\begin{figure}[!ht]
  \centering
  \includegraphics[page=5,scale=1.00]{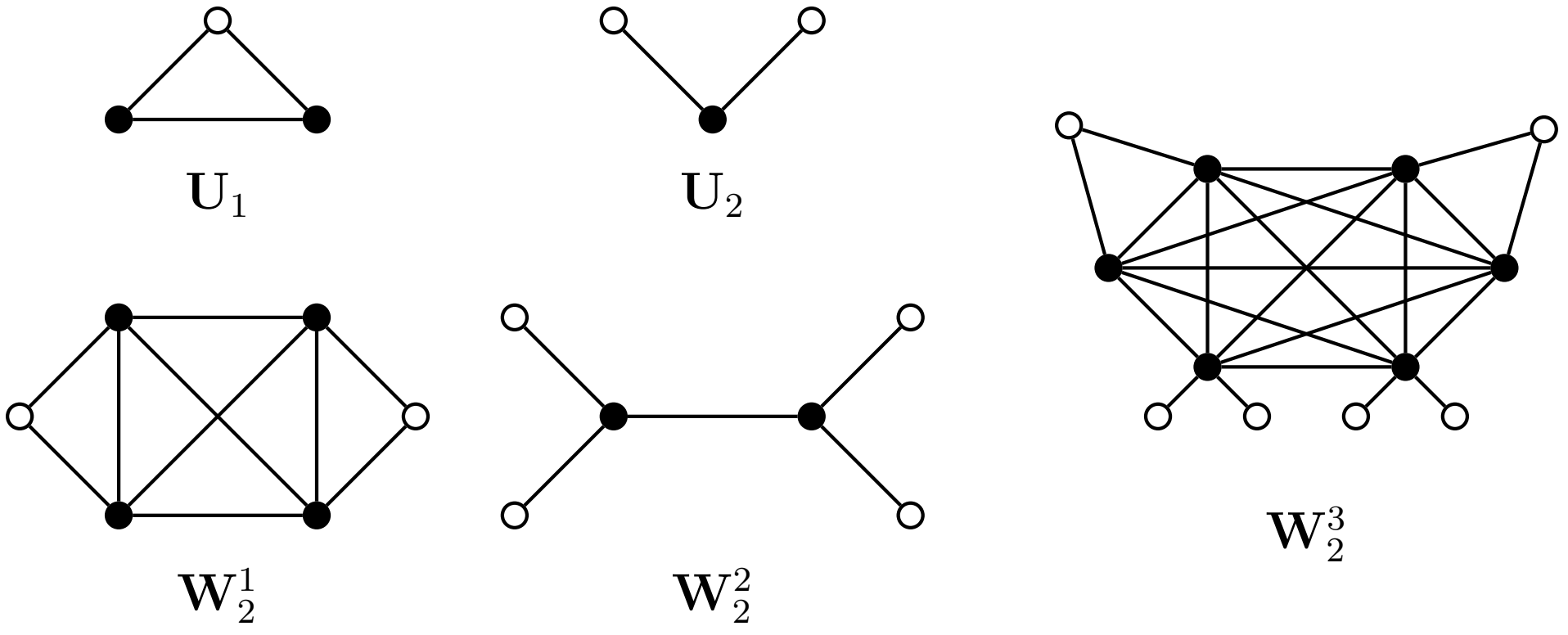}
  \caption{The marked graphs $\mathbf{U}_1$, $\mathbf{U}_2$, $\mathbf{W}_{2}^{1}$, $\mathbf{W}_{2}^{2}$ and $\mathbf{W}_{2}^3$.}
  \label{fig:7}
\end{figure}

\begin{example} Let $W_{n}^1$ denote the graph of order  $3n$ obtained from  $K_{2n}$ with vertex set $\{1,2,\ldots,2n\}$ by inserting $n$ isolated vertices, say $w_1,w_2,\ldots,w_n$, and then joining $w_i$ to the vertices $2i-1$ and $2i$.
	Let $W_{n}^2$ denote the graph of order  $3n$ obtained  by attaching two pendant vertices at each vertex of~$K_n$.

	In \cite{Liu2024} and \cite{Pet1998}, the inequalities $ \lambda_2(W_n^1) \leq 1$ and $ \lambda_2(W_n^2) \leq 1$ are verified by means of algebraic computation, and here we offer an alternative proof. Let $\mathbf{U}_1$ and $\mathbf{U}_2$ be the marked graphs depicted in Figure~\ref{fig:7}. It can be verified that $\mathbf{U}_1=\overline{\mathbf{L}^*(K_1(1))}$ and $\mathbf{U}_2=\overline{\mathbf{L}^*(K_2)}$. For $i\in\{1,2\}$, let $\mathbf{W}_{n}^{i}$ be the marked join of $n$ copies of $\mathbf{U}_i$. Note that $W_n^i$ is the underlying  graph of $\mathbf{W}_{n}^{i}$. By Corollary~\ref{cor2}, it follows that $\lambda_2(W_n^i) \leq 1$.
\end{example}

Moreover, by combining distinct marked graphs, one can further generate a structurally rich family of graphs satisfying $\lambda_2 \le 1$.
\begin{example}
Let $\mathbf{W}_{n}^3=\mathbf{W}_{n}^{1}\widehat{\vee}\mathbf{W}_{n}^{2}$, where $\mathbf{W}_{n}^{1}$ and $\mathbf{W}_{n}^{2}$ are taken on disjoint vertex sets. Equivalently, $\mathbf{W}_{n}^3$ is the marked join of $n$ copies of $\mathbf{U}_{1}$ and $n$ copies of $\mathbf{U}_{2}$. The underlying graph $W_n^3$ is of order $6n$ obtained from $W_n^1\cup W_n^2$ by joining every vertex of the $K_{2n}$-subgraph in $W_n^1$ to every vertex of the $K_n$-subgraph in $W_n^2$. Again, Corollary~\ref{cor2} implies that $\lambda_2(W_n^3) \leq 1$.	
\end{example}

\section*{Declaration of competing interest}

The authors declare that they have no known competing financial interests or personal relationships that could have
appeared to influence the work reported in this paper.

\section*{Data availability}
The computations presented in this study are implemented by using free, open-source mathematical software system SageMath~\cite{sagemath} and are publicly accessible in the form of Jupyter Notebook files. The corresponding repository is hosted on GitHub at: \url{https://github.com/Chaohui-Chen/CGSS_marked_graphs}.

\begin{thebibliography}{99}
	
	\bibitem{Bel2016}
	F. Belardo, T. Pisanski, and S.K. Simi\'c.
	\newblock On graphs whose least eigenvalue is greater than $-2$.
	\newblock \emph{Linear Multilinear Algebra} 64 (2016) 1570--1582.  
	
    \bibitem{Bou2002}
	N. Bourbaki.
	\newblock \emph{Lie Groups and Lie Algebras: Chapters 4--6}.
	\newblock  Elements of Mathematics, Springer-Verlag, Berlin, 2002.
    
    \bibitem{Weyl1912}
	A.E. Brouwer, and  W.H. Haemers.  	\newblock \emph{Spectra of Graphs}.	\newblock   Springer, New York, 2012.
	
	\bibitem{Cam1976}
	P.J. Cameron, J. M. Goethals, J.J. Seidel, and E.E. Shult.
	\newblock Line graphs, root systems, and elliptic geometry.
	\newblock \emph{J. Algebra} 43 (1976) 305--327.  
	
	\bibitem{Cas1997}
	J.W.S. Cassels.
	\newblock \emph{An Introduction to the Geometry of Numbers}.
	\newblock  Classics in
	Mathematics, Springer-Verlag, Berlin, 1997.
	
	\bibitem{Chen2025}
	C. Chen, H. Shan, M. Liu, and Z. Stani\'c.
	\newblock Connected claw-free graphs whose second largest eigenvalue does not exceed $1$.
	\emph{Discrete Math.} 349 (2026) 115079. 
	
	\bibitem{Cheng2018}
	X.-M. Cheng, G.R.W. Greaves, and J.H. Koolen.
	\newblock Graphs with three eigenvalues and second largest eigenvalue at most 1.
	\newblock \emph{J. Combin. Theory Ser. B} 129 (2018) 55--78.  
	
	\bibitem{Cheng2025}
	K. Cheng.
	\newblock On 
	$K_4$-minor free graphs whose second largest eigenvalue does not exceed 1.
	\newblock  \emph{European J. Combin.}  136 (2026), 104397.   
	
	\bibitem{Cvet1982}
	D. Cvetkovi\'c.
	\newblock On graphs whose second largest eigenvalue does not exceed 1.
	\newblock \emph{Publ. Inst. Math. $($Beograd$)$} 31 (1982) 15--20.  
		
	\bibitem{Cvet2004}
	D. Cvetkovi\'c, P. Rowlinson, and S. Simi\'c.
	\newblock \emph{Spectral Generalizations of Line Graphs: On Graphs with Least Eigenvalue $-2$}.
	\newblock  Cambridge University Press, Cambridge, 2004.
	
	\bibitem{Cvet2010}
	D. Cvetkovi\'c, P. Rowlinson, and S. Simi\'c.
	\newblock \emph{An Introduction to the Theory of Graph Spectra}.
	\newblock  Cambridge University Press, Cambridge, 2010.
	
	\bibitem{Cvet2015}
	D. Cvetkovi\'c, P. Rowlinson, and S. Simi\'c.
	\newblock Graphs with least eigenvalue $-2$: ten years on.
	\newblock \emph{Linear Algebra Appl.} 484 (2015) 504--539.  
	
	\bibitem{Geng2024}
	Z. Geng, J.-M. Guo, and Z. Wang.
	\newblock The $T_{1,2}$-free planar graphs whose second largest eigenvalue does not exceed 1.
	\newblock \emph{Discrete Math.} 347 (2024) 113999.  
	
	\bibitem{Grea2015}
	G. Greaves, J.  Koolen, A. Munemasa, Y. Sano, and T. Taniguchi.
	\newblock Edge-signed graphs with smallest eigenvalue greater than $-2$.
	\newblock \emph{J. Combin. Theory Ser. B} 110 (2015) 90--111.  
	
	\bibitem{Hal2015}
	B.C. Hall.
	\newblock \emph{Lie Groups, Lie Algebras, and Representations: An Elementary Introduction}.
	\newblock Grad. Texts Math. vol. 222. Springer, Cham, 2nd ed., 2015.  
	
	\bibitem{Hof1968}
	A.J. Hoffman.
	\newblock Some recent results on spectral properties of graphs.
	\newblock In \emph{Beitr. Graphentheorie, Int. Kolloquium Manebach} (DDR) 1967, pp. 75--80. 1968.
	
	\bibitem{Hum1972}
	J.E. Humphreys.
	\newblock \emph{Introduction to Lie Algebras and Representation Theory}.
	\newblock Grad. Texts Math. vol. 9. Springer, Cham, 1972.  
	
	\bibitem{Koo2021}
	J.H. Koolen, M.-Y. Cao, and Q. Yang.
	\newblock Recent progress on graphs with fixed smallest adjacency eigenvalue: a survey.
	\newblock \emph{Graphs Combin.} 37 (2021) 1139--1178.  
	
	\bibitem{Li2023}
	S. Li and W. Sun.
	\newblock Characterization of outerplanar graphs whose second largest eigenvalue is at most 1.
	\newblock \emph{Result. Math.} 78 (2023) 104.  
	
	\bibitem{Liu2021}
	M. Liu, C. Chen, and Z. Stani\'c.
	\newblock On graphs whose second largest eigenvalue is at most 1.
	\newblock \emph{European J. Combin.} 97 (2021) 103385.  
	
	\bibitem{Liu2024}
	M. Liu, C. Chen, and Z. Stani\'c.
	\newblock Connected $(K_4-e)$-free graphs whose second largest eigenvalue does not exceed 1.
	\newblock \emph{European J. Combin.} 115 (2024) 103775.  
	
	\bibitem{Liu2025}
	M. Liu, C. Chen, Z. Stani\'c, and H. Shan.
	\newblock Graphs with large clique number whose second largest eigenvalue does not exceed $(\sqrt{5}-1)/2$.
	\newblock \emph{Electron. J. Combin.} 32 (2025) \#P2.13.  
	
	\bibitem{Pet1991}
	M. Petrovi\'c.
	\newblock On graphs with exactly one eigenvalue less than $-1$.
	\newblock \emph{J. Combin. Theory Ser. B} 52 (1991) 102--112.  
	
	\bibitem{Pet1998}
	M. Petrovi\'c and B. Mileki\'c.
	\newblock On the second largest eigenvalue of line graphs.
	\newblock \emph{J. Graph Theory} 27 (1998) 61--66.  
	
	\bibitem{Pet2000}
	M. Petrovi\'c and B. Mileki\'c.
	\newblock Generalized line graphs with the second largest eigenvalue at most 1.
	\newblock \emph{Publ. Inst. Math. $($Beograd$)$} 68 (2000) 37--45.
	 
	
	\bibitem{Row2022}
	P. Rowlinson and Z. Stani\'c.
	\newblock Signed graphs whose spectrum is bounded by $-2$.
	\newblock \emph{Appl. Math. Comput.} 423 (2022) 126991.  
	
	\bibitem{sagemath}
	The Sage Developers.
	\newblock \emph{SageMath, The Sage Mathematics Software System}, version 10.7 [software], 2025.
	 

	
	\bibitem{Stanic2007}
	Z. Stani\'c.
	\newblock On graphs whose second largest eigenvalue equals 1--the star complement technique.
	\newblock \emph{Linear Algebra Appl.} 420 (2007) 700--710.  
	
	\bibitem{Stanic2009}
	Z. Stani\'c.
	\newblock On nested split graphs whose second largest eigenvalue is less than 1.
	\newblock \emph{Linear Algebra Appl.} 430 (2009) 2200--2211.  
	
	\bibitem{Stanic2010}
	Z. Stani\'c.
	\newblock On regular graphs and coronas whose second largest eigenvalue does not exceed 1.
	\newblock \emph{Linear Multilinear Algebra} 58 (2010) 545--554.  
	
	\bibitem{ifge}
	Z. Stani\'c.
	\newblock \emph{Inequalities for Graph Eigenvalues}.
	\newblock Cambridge University Press, 2015.
	

\end{thebibliography}

\end{document}